\documentclass[aos,MSNbibl,nameyear,dvips]{arximspdf}

\doi{10.1214/12-AOS971} 
\volume{40}
\issue{1}
\pubyear{2012}
\firstpage{561}
\lastpage{592}

\makeatletter

\newcommand{\cal}{\mathcal}

\newcommand{\R}{\mathbb{R}}
\newcommand{\reals}{\R}

\newcommand{\bQ}{{\mathbf Q}}
\newcommand{\bp}{{\mathbf p}}
\newcommand{\bu}{{\mathbf u}}
\newcommand{\D}{D}
\newcommand{\E}{{\mathrm{E}}}

\newcommand{\defeq}{:=}
\newcommand{\half}{\frac{1}{2}}
\newcommand{\cd}{|}
\newcommand{\bldq}{\mathbf{q}}

\newtheorem{theorem}{Theorem}[section]
\newtheorem{lem}[theorem]{Lemma}

\newproclaim{expl}{Example}[section]
\newproclaim{definer}{Definition}[section]
\newproclaim{cond}{Condition}[section]

\newtheorem{cor}[theorem]{Corollary}

\makeatother

\begin{document}
\begin{frontmatter}

\title{Proper local scoring rules}
\runtitle{Proper local scoring rules}

\begin{aug}
\author[A]{\fnms{Matthew} \snm{Parry}\thanksref{t1}\ead[label=e1]{mparry@maths.otago.ac.nz}},
\author[B]{\fnms{A. Philip} \snm{Dawid}\corref{}\ead[label=e2]{A.P.Dawid@statslab.cam.ac.uk}}
\and
\author[C]{\fnms{Steffen} \snm{Lauritzen}\ead[label=e3]{steffen@stats.ox.ac.uk}\ead[label=u1,url]{http://www.foo.com}}
\runauthor{M. Parry, A. P. Dawid and S. Lauritzen}
\affiliation{University of Otago, University of Cambridge and
University of Oxford}
\address[A]{M. Parry\\
Department of Mathematics\\
\quad and Statistics\\
University of Otago\\
P.O. Box 56\\
Dunedin 9054\\
New Zealand\\
\printead{e1}}
\address[B]{A. P. Dawid\\
Statistical Laboratory\\
Centre for Mathematical Sciences\\
University of Cambridge \\
Wilberforce Road, Cambridge CB3 0WB\\
United Kingdom\\
\printead{e2}}
\address[C]{S. Lauritzen\\
Department of Statistics \\
University of Oxford\\
South Parks Road\\
Oxford OX1 3TG\\
United Kingdom\\
\printead{e3}\\
\printead{u1}} 
\end{aug}

\thankstext{t1}{Supported by EPSRC Statistics Mobility
Fellowship EP/E009670.}

\received{\smonth{1} \syear{2011}}
\revised{\smonth{1} \syear{2012}}

%
\begin{abstract}
We investigate proper scoring rules for continuous distributions on
the real line. It is known that the log score is the only such rule
that depends on the quoted density only through its value at the
outcome that materializes. Here we allow further dependence on a~finite number $m$ of derivatives of the density at the outcome, and
describe a~large class of such \textit{$m$-local} proper scoring
rules: these exist for all even $m$ but no odd $m$. We further show
that for $m\geq2$ all such $m$-local rules can be computed without
knowledge of the normalizing constant of the distribution.
\end{abstract}

%

\end{frontmatter}

\section{Introduction}
\label{secintro}
A~\textit{scoring rule} $S(x, Q)$ is a~loss function measuring the
quality of a~quoted distribution $Q$, for an uncertain quantity $X$,
when the realized value of $X$ is $x$. It is \textit{proper} if it
encourages honesty in the sense that the expected score $\E_{X \sim
P}S(X,Q)$, where $X$ has distribution $P$, is minimized by the
choice $Q=P$.

Traditionally, a~scoring rule has been termed \textit{local} if it
depends on the density function $q(\cdot)$ of $Q$ only through its
value, $q(x)$, at $x$. With this definition, any proper local scoring
rule is equivalent to the \textit{log score}, $S(x,Q) = -\ln q(x)$.
However, we can weaken the locality condition by allowing further
dependence on a~finite number $m$ of derivatives of $q(\cdot)$ at $x$,
and this introduces many further possibilities. We term $m$ the
\textit{order} of the rule.

In this paper we describe a~large class of such \textit{order-$m$ proper
local scoring rules} for densities on the real line. These turn
out to depend on the density~$q(\cdot)$ in a~way that is insensitive
to a~multiplicative constant, and hence can be computed without
knowledge of the normalizing constant of $q$.

\citet{Hyvarinen2005es} proposed a~method for approximating a~distribution $P$ on $\mathcal{X}=\R^k$ by a~distribution $Q$ in a~specified family $\mathcal{P}$ of distributions by minimizing $d(P,Q)$
over $Q\in\mathcal{P}$, where
%
\begin{equation}
\label{eqhyvd}
d(P,Q)=\frac{1}{2}\int {\mathrm{d}}x\, p(x) |\nabla\ln
p(x)-\nabla\ln
q(x) |^{2}
\end{equation}
with $\nabla$ denoting gradient. Since $q$ enters this expression only
through $\nabla\ln q$, it is clear that the minimization only
requires knowledge of $q$ up to a~multiplicative factor. Using
integration by parts, \citet{Hyvarinen2005es} further showed that
minimization of the divergence $d(P,Q)$ in (\ref{eqhyvd}) is
equivalent to mimimizing
\[
S(P,Q)=\E_P\bigl\{\Delta\ln q(X)+\tfrac{1}{2}|\nabla\ln
q(X)|^{2}\bigr\}
\]
[where $\Delta$ denotes the Laplacian operator $\sum_{i=1}^k
\partial^2/(\partial x_i)^2$], which is a~scoring rule of the type
discussed in this paper: see Section~\ref{sechyvscore} below.

The plan of the paper is as follows. In Section~\ref{secsr} we introduce
proper scoring rules, with some examples and applications.
Section~\ref{seclocalscoringrules} formalizes the notion of a~local
function, its representations and derivatives. In Section~\ref{secvaranal} we
apply integration by parts and the calculus of variations to develop a~``key equation,'' which is further investigated in Section~\ref{secsolve}
through an analysis of fundamental differential operators associated
with local functions. Section~\ref{sechomscore} describes the solutions to
the key equation, which we term ``key local scoring rules,'' in terms
of a~homogeneous function $\phi$. In Section~\ref{secnonunique} we point out
that distinct choices of $\phi$ can generate the same scoring rule,
and consider some implications; in particular, we show that key
$m$-local scoring rules exist for any even order $m$, but for no odd
order. Section~\ref{secprop} examines when this construction does indeed
yield a~proper local scoring rule, concavity of $\phi$ being crucial.
Section~\ref{secboundary} devotes further attention to boundary terms
arising in the integration by parts. In Section~\ref{sectransform} we study
how the problem and its solution transform under an invertible mapping
of the sample space, and develop an invariant formulation.

\subsection{Related work}
\label{secrelated}

In this paper we are concerned with characterizing $m$-local proper
scoring rules, for all orders $m$. Since there are no such rules of
order 1, order~2 scoring rules constitute the simplest nontrivial
case, and as such are likely to be the most useful in practice. In
a~companion paper to this one, \citet{ehmgneitingorder2} conduct a~deep
investigation of order~2 proper local rules, using an elegant
construction complementary to ours. They also describe a~general
class of densities for which the boundary terms vanish.

The present paper confines attention to absolutely continuous
distributions on the real line. The notion of local scoring rule has
an interesting analogue for a~discrete sample space equipped with a~given neighborhood structure. The theory for that case is developed
in an accompanying paper [\citet{apdsllmfpdiscrete}]; it exhibits
both close parallels with, and important differences from, the
continuous case considered here.

\section{Scoring rules}
\label{secsr}
Suppose You are required to express Your uncertainty about an
unobserved quantity $X\in{\mathcal X}$ by quoting a~distribution~$Q$
over ${\mathcal X}$, after which Nature will reveal the value $x$ of~$X$.
A~\textit{scoring rule} or \textit{score}~$S$
[\citet{apdencprobfore}] is a~special kind of loss function, intended
to measure the quality of your quote $Q$ in the light of the realized
outcome $x$: $S(x,Q)$ is a~real number interpreted as the loss You
will suffer in this case. The principles of Bayesian decision theory
[\citet{savagebook}] now enjoin You to minimize Your expected loss. If
Your actual beliefs about $X$ are described by a~probability
distribution $P$, You should thus quote that $Q$ that minimizes
$S(P,Q) \defeq\E_{X \sim P}{S(X,Q)}$. The scoring rule $S$ is termed
\textit{proper} (relative to a~class ${\mathcal P}$ of distributions
over ${\mathcal X}$) when, for any fixed $P\in{\mathcal P}$, the
minimum over $Q\in{\cal P}$ is achieved at $Q=P$; it is \textit{strictly
proper} when, further, this minimum is unique. Thus, under a~proper scoring rule, honesty is the best policy.

Associated with any proper scoring rule $S$ are a~(\textit{generalized})
\textit{entropy function $H(P) \defeq S(P,P)$} and a~\textit{divergence
function} $d(P,Q) \defeq S(P,Q)-H(P)$. Under suitable technical
conditions, proper scoring rules and their associated entropy
functions and divergence functions enjoy certain properties that serve
to characterize such ``coherent'' constructions [\citet{apd94}]:
$S(P,Q)$ is affine in~$P$ and is minimized in $Q$ at $Q=P$; $H(P)$ is
concave in~$P$; $d(P,Q)-d(P,Q_{0})$ is affine in $P$, and $d(P,Q)\geq
0$, with equality achieved at $Q=P$.

If two scoring rules differ by a~function of $x$ only, then they will
yield the identical divergence function. In this case we will term
them \textit{equivalent} [note that this is a~more specialized usage
than that of \citet{apd94}].

A~fairly arbitrary statistical decision problem can be reduced to one
based on a~proper scoring rule. Let
$L\dvtx\mathcal{X}\times\mathcal{A}\rightarrow\R$ be a~loss function,
defined for outcome space $\mathcal{X}$ and action space
$\mathcal{A}$. Letting $\mathcal{P}$ be a~class of distributions over
$\mathcal{X}$ such that $L(P,a) \defeq\mbox{E}_{X\sim P}L(X,a)$ exists
for all $a\in\mathcal{A}$ and $P\in\mathcal{P}$, define, for
$P,Q\in\mathcal{P}$ and $x\in\mathcal{X}$,
%
\begin{equation}
S(x,Q)\defeq L(x,a_{Q}),
\end{equation}
where $a_{P} \defeq\arg\inf_{a\in\mathcal{A}}L(P,a)$ is a~Bayes act
with respect to $P$ (supposed to exist, and selected arbitrarily if
nonunique). Then $S$ is readily seen to be a~proper scoring rule,
and the associated entropy function is just the Bayes loss: $H(P) =
\inf_{a\in{\mathcal A}}L(P,a)$.

In this paper we focus attention on the case that $\mathcal{X}$ is an
interval on the real line and any $Q \in{\mathcal P}$ has a~density
$q(\cdot)$ with respect to Lebesgue measure on $\mathcal{X}$. We may
then define $S(x,Q)$ in terms of $q$. However, since $q$ is only
defined almost everywhere we must take care that any manipulations
performed either involve a~preferred version of $q$, or yield the same
answer when $q$ is changed on a~null set. This will always be the
case in this paper.

\subsection{Bregman scoring rule}

Since any decision problem generates a~proper scoring rule there is a~very great number of these. Certain forms are of special interest or
simplicity. Here we describe one important class of such rules for
the case that every $Q \in\mathcal{P}$ has a~density function
$q(\cdot)$ with respect to a~dominating measure $\mu$ over
$\mathcal{X}$.

Let $\phi\dvtx\R^{+}\rightarrow\R$ be concave and differentiable. The
associated (\textit{separable}) \textit{Bregman scoring rule} is defined by
%
\begin{equation}
S(x,Q) \defeq
\phi'\{q(x)\}+\int {\mathrm{d}}\mu(y)\,[\phi\{q(y)\}-q(y)\phi'\{
q(y)\}].
\end{equation}
It can be shown that these are the only proper scoring rules having
the form $S(x,Q) = \xi\{q(x)\} - k(Q)$ [\citet{apdpsr}].

Taking expectations, we obtain
\[ S(P,Q)=\int {\mathrm{d}}\mu(x)\,[
\{p(x)-q(x)\}\phi'\{q(x)\}+\phi\{q(x)\}].
\]
It follows that $H(P)=\int {\mathrm{d}}\mu(x) \,\phi\{p(x)\}$ and so, assuming
$H(P)$ is finite, the corresponding (separable) \textit{Bregman
divergence} [\citet{bregman67}, \citet{csiszar91}]---also termed
\textit{$U$-divergence} [\citet{eguchi2008}]---is
%
\begin{equation}
\label{eqbregd}
d(P,Q)=\int {\mathrm{d}}\mu(x)\,
\bigl([\phi\{q(x)\}+\{p(x)-q(x)\}\phi'\{q(x)\}]-\phi\{p(x)\}\bigr).
\end{equation}
The integrand is nonnegative by concavity of $\phi$. Therefore, the
separable Bregman scoring rule is a~proper scoring rule, and strictly
proper if $\phi$ is strictly concave.

\subsection{Extended Bregman score}
\label{secextbreg}

A~straightforward generalization of the above Bregman construction is
obtained on replacing $\phi\dvtx\mathbb{R}^+\rightarrow\mathbb{R}$
throughout by $\phi\dvtx\mathcal{X}\times
\mathbb{R}^+\rightarrow\mathbb{R}$, such that, for each
$x\in\mathcal{X}$, $\phi(x, \cdot)\dvtx\mathbb{R}^+\rightarrow\mathbb{R}$
is concave. Such \textit{extended Bregman} rules are the only proper
scoring rules of the form $S(x,Q) = \xi\{x,q(x)\} - k(Q)$
[\citet{apdpsr}].

\subsection{Log score}
\label{seclogscore}

For $\phi(s) \equiv- s\ln s$ we obtain the \textit{logarithmic scoring
rule}, or \textit{log score}, defined by
\[
S(x,Q)=-\ln q(x).
\]
This is essentially the only scoring rule of the form $S(x,Q) =
\xi\{x,q(x)\}$ [Bernardo (\citeyear{Bernardo1979ex}), \citet{apdpsr}]. For this case we
obtain
\[
H(P)=-\int {\mathrm{d}}\mu(x)\, p(x)\ln p(x),
\]
the \textit{Shannon entropy}, and
%
\begin{equation}
d(P,Q)=\int {\mathrm{d}}\mu(x)\,
p(x)\ln\frac{p(x)}{q(x)},
\end{equation}
the \textit{Kullback--Leibler divergence}.

\subsection{Parameter estimation}
\label{secparest}

Let $\mathcal{Q}=\{Q_{\theta}\}\subseteq\mathcal{P}$ be a~smooth
parametric family of distributions. Given data $(x_1, \ldots, x_n)$
in $\mathcal{X}$ with empirical distribution~$\widehat{P}$, one way to
estimate $\theta$ is by minimizing some divergence criterion:
$\widehat{\theta}\defeq\arg\min_{\theta}d(\widehat{P},Q_{\theta})$.
When the divergence function is derived from a~scoring rule, this is
equivalent to minimizing the total empirical score:
%
\begin{equation}
\widehat{\theta}=\mathop{\arg\min}_{\theta}\sum_{i=1}^n S(x_i,Q_{\theta}),
\end{equation}
in which form it remains meaningful even if
$\widehat{P}\notin\mathcal{P}$, when $d(\widehat{P},Q_{\theta})$ is
undefined. The corresponding estimating equation is
%
\begin{equation}
\label{eqestimating}
\sum_{i=1}^n \sigma(x_i,\theta)=0
\end{equation}
with $\sigma(x, \theta) \defeq{\partial
S(x,Q_\theta)}/{\partial\theta}$. For a~proper scoring rule it is
straightforward to show that the estimating equation
(\ref{eqestimating}) is unbiased [\citet{Dawid2005ge}] and, as a~result, $\widehat{\theta}$ is typically consistent, though not
necessarily efficient; it may also display some degree of robustness.
Equation (\ref{eqestimating}) delivers an \textit{$M$-estimator}
[\citet{Huber1981}, \citet{HampelETAL1986}]. Statistical properties of the
estimator are considered by \citet{eguchi2008} for the special case
of minimum Bregman ($U$-) divergence estimation, and readily extend to
more general cases.

\subsection{Hyv\"arinen scoring rule}
\label{sechyvscore}

Hyv{\"a}rinen (\citeyear{Hyvarinen2005es}) showed that minimization of the divergence
$d(P,Q)$ in (\ref{eqhyvd}) is equivalent to minimizing the empirical
score for the scoring rule
%
\begin{equation}
\label{eqgenhy}
S(x,Q)=\Delta\ln q(x)+\tfrac{1}{2}|\nabla\ln
q(x)|^{2}.
\end{equation}
This is valid in the case where $\mathcal{X}=\mathbb{R}^{k}$ and
${\cal P}$ consists of distributions $P$ whose Lebesgue density
$p(\cdot)$ is a~twice continuously differentiable function of $x$
satisfying $\nabla\ln p\rightarrow0$ as $|x|\rightarrow\infty$. For
$k=1$ we get
%
\begin{equation}
\label{eqhyv1}
S(x,Q) = \frac{q''(x)}{q(x)} - \half\biggl\{\frac{q'(x)}{q(x)}
\biggr\}^2.
\end{equation}

\citet{Dawid2005ge} showed that, with some reinterpretation, the
formula (\ref{eqgenhy}) defines a~proper scoring rule in the more
general case of an outcome space $\mathcal{X}$ that is a~Riemannian
manifold. Now $q(\cdot)$ denotes the \textit{natural density} $\mathrm{d}
Q/\mathrm{d}\mu$ of $Q$ with respect to the associated volume measure
$\mu$ on $\mathcal{X}$; $\nabla$ denotes \textit{natural gradient};
$\Delta$ is the \textit{Laplace--Beltrami operator}; and
$|u|^{2}=\langle u,u \rangle$ is the squared norm defined by the
metric tensor. We impose the restriction $P,Q \in\mathcal{P}$, where
$P\in{\cal P}$ if $\nabla\ln p(x)\rightarrow0$ as $x$ approaches the
boundary of $\mathcal{X}$.

On applying Stokes's theorem (again essentially integration by parts)
and noting that boundary terms vanish, we can express the expected
score as
\[
S(P,Q)=\frac{1}{2}\int {\mathrm{d}}\mu(x)\, p(x) \langle\nabla\ln
q(x)-2\nabla\ln p(x), \nabla\ln q(x)\rangle.
\]
The entropy is thus
$H(P)=-\frac{1}{2}\int {\mathrm{d}}\mu(x)\,p(x) |\nabla\ln
p(x)|^{2}$ and so the associated divergence is essentially
that used by Hyv{\"a}rinen:
%
\begin{equation}
\label{eqgenhyvd}
d(P,Q)=\frac{1}{2}\int {\mathrm{d}}\mu(x)\,p(x) |\nabla\ln
p(x)-\nabla\ln
q(x) |^{2},
\end{equation}
which is nonnegative and vanishes only when $Q=P$. It follows that
the scoring rule is strictly proper.

Although this scoring rule is not local in the strict sense, it
depends on $(x,Q)$ only through the first and second derivatives of
the density function~$q(\cdot)$ at the point $x$; it is \textit{local of
order $2$}, or \textit{$2$-local}, as defined below in
Section~\ref{seclocalscoringrules}.

Note that one does not need to know the volume measure $\mu$ to
calculate the divergence; formula (\ref{eqgenhyvd}) for $d(P,Q)$
yields the same result if we take $\mu$ to be any fixed underlying
measure, and interpret $p$ and $q$ as densities with respect to this.

\subsection{Homogeneity}
\label{sechyvhom}

An interesting and practically valuable property of the generalized
Hyv\"arinen scoring rule is that $S(x,Q)$ given by (\ref{eqgenhy}) is
\textit{homogeneous} in the density function $q(\cdot)$: it is formally
unchanged if $q(\cdot)$ is multiplied by a~positive constant, and so
can be computed even if we only know the density function up to a~scale factor. In particular, use of the estimating equation
(\ref{eqestimating}) does not require knowledge of the normalizing
constant (which is often hard to obtain) for densities in
$\mathcal{Q}$.
%
\begin{expl}
\label{exselection}
Consider the natural exponential family:
\[
q(x \cd\theta) = Z(\theta)^{-1} \exp\{a(x) + \theta x\}.
\]
Using the scoring rule (\ref{eqhyv1}) we obtain $S(x, Q_\theta) =
a''(x) + \half\{a'(x)+\theta\}^2$, so that $\sigma(x,\theta) = a'(x)
+ \theta$, and (\ref{eqestimating}) delivers the unbiased estimator
\[
\widehat\theta= - \sum_{i=1}^n a'(X_i)/n,
\]
which can be computed without knowledge of $Z(\theta)$. See also
Section 4 of \citet{Hyvarinenext}, where exponential families are
discussed.

Alternatively, we can work directly with the sufficient statistic $T
\defeq\sum_{i=1}^n X_i$, which has density of the form
\[
q_T(t \cd\theta) = Z(\theta)^{-n} \exp\{\alpha_n(t) + \theta t\}.
\]
Applying the above method to $q_T$ leads to the unbiased estimator
\mbox{$\widetilde\theta= -\alpha_n'(T)$}. This is the
\textit{maximum plausibility estimator} of $\theta$
[\citet{barndorff1976}]. Basing the estimate on the sufficient
statistic is more satisfying and better behaved from a~principled
point of view, but does require computation of the function
$\alpha_n(t)$, which involves an $n$-fold convolution of $e^{a(x)}$.

As an application, suppose $Q_\theta$ is obtained from the normal
distribution~$N(\theta,1)$ by retaining its outcome~$x$ with
probability~$k(x)$. We assume that~$k(x)$ is everywhere positive
and twice differentiable. The density is thus
\[
q(x \cd\theta) = \frac{k(x) \exp-\{(x-\theta)^2/2\}}
{\int k(y) \exp-\{(y-\theta)^2/2\} \,\mathrm{d}y}.
\]
Because of the complex dependence of the denominator on $\theta$,
the maximum likelihood estimate typically cannot be expressed in
closed form. However, using scoring rule (\ref{eqhyv1}) yields the
explicit unbiased estimator
\[
\widehat\theta= \sum_{i=1}^n\{x_i - \kappa'(x_i)\}/n
\]
with $\kappa(x) \defeq\ln k(x)$.
\end{expl}

The homogeneity property will be a~feature of all the new proper local
scoring rules we introduce here: see Section~\ref{sechomscore}.

\section{Local scoring rules}
\label{seclocalscoringrules}

We observed in Section~\ref{seclogscore} that the log score $S(x, Q) = -\ln
q(x)$ is essentially the only proper scoring rule that is
\textit{local}, that is, involves the density function $q(\cdot)$ of $Q$ only
through its value, $q(x)$, at the actually realized value $x$ of $X$.

We can, however, weaken the locality requirement, for example, by
allowing $S$ to depend on the values of $q(\cdot)$ in an infinitesimal
neighborhood of $x$. In this paper we describe a~class of scoring rules
that depend on the function $q(\cdot)$ only through its value and the
values of a~finite number $m$ of its derivatives at the
point~$x$---a~property we will refer to as \textit{locality of order~$m$}, or \textit{$m$-locality}.

We confine ourselves to the case that $\mathcal{X}$ is an open
interval in $\mathbb{R}$, possibly infinite or semi-infinite, and
$\mathcal{P}$ is a~class of distributions $Q$ over $\mathcal{X}$
having strictly positive Lebesgue density $q(\cdot)$ that is $m$-times
continuously differentiable.

\subsection{Local functions and scoring rules}
To study the properties of local scoring rules we need a~formal
definition of a~local function.
%
\begin{definer}
\label{deforder}
A~function $F\dvtx\mathcal{X}\times\mathcal{P}\to\R$ is said to be
\textit{local of order~$m$}, or~\textit{$m$-local}, if it can be
expressed in the form
\[
F(x,Q)=f\bigl\{x,q(x),q'(x),q''(x),\ldots,q^{(m)}(x)\bigr\},
\]
where $f\dvtx\mathcal{X}\times\bQ_m\to\R$, with $ \bQ_m \defeq\R^+
\times\R^{m}$, is a~real-valued infinitely differentiable function,
$q(\cdot)$ is the density function of $Q$, and a~prime ($'$) denotes
differentiation with respect to $x$. It is \textit{local} if it is
local of some finite order.
\end{definer}

We shall refer to such a~function $f$ as a~\textit{$q$-function}, and
say it is of order~$m$. When we do not need to specify the order $m$
of a~$q$-function $f$ we may write $f(x,\bldq )$ ($x\in
\mathcal{X}$, $\bldq
\in\bQ$), understanding $\bQ= \bQ_m$, $\bldq =
(q_0,\ldots,q_m)$.

A~scoring rule $S(x, Q)$ is \textit{$m$-local} if
%
\begin{equation}
\label{eqmscore}
S(x,Q)=s\bigl\{x,q(x),q'(x),q''(x),\ldots,q^{(m)}(x)\bigr\},
\end{equation}
where $s$ is a~$q$-function as above, so that it depends on the quoted
distribution $Q$ for $X$ only through the value and derivatives up to
order $m$ of the density $q(\cdot)$ of~$Q$, evaluated at the observed
value $x$ of $X$. The function $s$ is the \textit{score function}
of~$S$.

\subsection{Differentiation of local functions}
\label{secdiffop}
For a~local scoring rule $S$ given by (\ref{eqmscore}) we write
\[
S_{[j]}(x,Q) \defeq s_{[j]}\bigl\{x,q(x),q'(x),q''(x), \ldots,q^{(m)}(x)\bigr\},
\]
where $s_{[j]} \defeq\partial s/\partial q_j$, and similarly
\[
S_{[x]}(x,Q) \defeq s_{[x]}\bigl\{x, q(x),\ldots, q^{(m)}(x)\bigr\},
\]
where $s_{[x]} \defeq\partial s/\partial x$. Then if ${\mathrm{d}}S/{\mathrm{d}}x$
denotes the derivative of $S(x,Q)$ with respect to $x$ for fixed $Q$,
we have
%
\begin{equation}
\label{eqtotal}
\frac{{\mathrm{d}}S}{{\mathrm{d}}x} = S_{[x]}(x,Q) +
\sum_{j \geq0}
q^{(j+1)}(x)   S_{[j]}(x,Q).
\end{equation}
For $S$ of order $m$, the series in (\ref{eqtotal}) terminates at $j=m$.

Motivated by (\ref{eqtotal}), we introduce a~linear differential
operator $\D$ acting on $q$-functions by
%
\begin{equation}
\label{eqtotally}
\D\defeq\frac{\partial}{\partial x} + \sum_{j\geq0} q_{j+1}
\,\frac{\partial}{\partial q_j}.
\end{equation}
For $f$ of order $m$, the series for $\D f$ obtained from
(\ref{eqtotally}) terminates at $j=m$, and $\D f$ is then of order
$m+1$.

The operator $\D$ thus represents the total derivative of the local
function for fixed~$Q$:
\[
\frac{{\mathrm{d}}S}{{\mathrm{d}}x} = (\D s)\bigl\{x, q(x), \ldots,
q^{(m+1)}(x)\bigr\},
\]
where $s$ is the score function of $S$.

In the light of the interpretation of $D$ as $\mathrm{d}/\mathrm{d}x$, the following
result is unsurprising:
%
\begin{lem}
\label{lemdint}
For a~$q$-function $f$, $Df = 0$ if and only if $f$ is
constant.
\end{lem}
\begin{pf}
``If'' is trivial. For ``only if,'' suppose $f$ is of order $\leq
m$. The only term in $Df$ involving $q_{m+1}$ is $q_{m+1} f_{[m]}$,
so that $Df = 0 \Rightarrow f_{[m]} = 0$, whence $f$ is of order at
most $m-1$. Repeating this argument, $f$ must be of order 0, that is,
of the form $f(x)$. Then $0=Df = f'(x)$, so finally $f$ must be a~constant.
\end{pf}

\section{Variational analysis}
\label{secvaranal}
We are interested in constructing proper local scoring rules. Ideally
we would develop sufficient conditions on the score function $s$ and
the family $\mathcal{P}$ to ensure that, for any $P\in\mathcal{P}$,
$S(P,Q)=\int {\mathrm{d}}x\, p(x)  s\{x,q(x)$, $q'(x),q''(x),
\ldots,q^{(m)}(x)\}$ is minimized, over $Q \in\mathcal{P}$, at $Q=P$.
Initially, however, we shall merely develop, in a~somewhat heuristic
fashion, conditions sufficient to ensure that, for all $P\in{\cal P}$,
$Q=P$ will be a~\textit{stationary point} of $S(P,Q)$---a~property we
shall describe by saying that~$S$ is a~\textit{weakly proper} scoring
rule. Given any $S$ satisfying these conditions, further attention
will be required to check whether or not it is in fact proper; this
will be taken up in Section~\ref{secprop} below.

To address this problem we adopt the methods of variational calculus
[\citet{Troutman1983}, \citet{vanbrunt04}]. Suppose that, at $Q=P$, $S(P,Q)$
is stationary under an arbitrary infinitesimal variation $\delta
q(\cdot)$ of $q(\cdot)$, subject to the requirement that $q(\cdot) +
\delta q(\cdot)$ be a~probability density. That is,
%
\begin{eqnarray}
\label{eqvariation}\quad
&&\delta\biggl\{\int {\mathrm{d}}x\, p(x)
s\bigl\{x,q(x),q'(x),q''(x),\ldots,q^{(m)}(x)\bigr\}
+\lambda_P \int {\mathrm{d}}x\, q(x)\biggr\}\bigg|_{q=p}
\nonumber\\[-8pt]\\[-8pt]
&&\qquad = 0,
\nonumber
\end{eqnarray}
where $\lambda_P$ is a~Lagrange multiplier for the normalization
constraint $\int{\mathrm{d}}x\,q(x) = 1$. The left-hand side of
(\ref{eqvariation}), evaluated with $P = Q$, is
%
\begin{equation}
\label{eqvarn}
\int {\mathrm{d}}x\,\Biggl\{\sum_{k=0}^{m} \delta
q^{(k)}(x)  q(x) S_{[k]}(x,Q) + \lambda_Q  \delta q(x)\Biggr\},
\end{equation}
and this is to vanish for arbitrary infinitesimal $\delta q(\cdot)$
and suitable $\lambda_Q$.

We evaluate the integral of the $k$th term of the sum in (\ref{eqvarn})
using the general formula for repeated integration by parts:
%
\begin{eqnarray}
\label{eqfullparts}
&&(-1)^k \int_-^+ {\mathrm{d}}x\, FG^{(k)}\nonumber\\[-8pt]\\[-8pt]
&&\qquad =\int_-^+ {\mathrm{d}}x\,
G F^{(k)} -
 \sum_{r=0}^{k-1}(-1)^{k-1-r}\bigl\{G^{(k-1-r)} F^{(r)}
\bigr\}
\big|_-^+,\nonumber
\end{eqnarray}
where $F^{(k)}$ denotes the $k$th derivative of $F$ with respect to
$x$. The first term on the right-hand side of\vadjust{\goodbreak} (\ref{eqfullparts}) is
the \textit{integral term}; the remaining terms are \textit{boundary
terms}, these being evaluated, if necessary, as limits as we
approach, from within, the end-points (denoted by $-$ and $+$) of the
interval ${\cal X} \subseteq\R$.

Setting $G = q(x) S_{[k]}(x, Q)$, $F=\delta q(x)$, we obtain
%
\begin{eqnarray}
\label{eqfullvariationq}
&&\int_-^+  {\mathrm{d}}x\, q(x) S_{[k]}(x,Q)\delta
q^{(k)}(x) \nonumber\\
&&\qquad=\int_-^+ {{\mathrm{d}}x}\, (-1)^k \delta q(x)\,
\frac{{\mathrm{d}}^{k}}{{\mathrm{d}}x^{k}}\,\bigl\{q(x)
S_{[k]}(x,Q)\bigr\}\\
&&\qquad\quad{}+ \sum_{r=0}^{k-1} (-1)^{k-1-r}\,
\frac{{\mathrm{d}}^{k-1-r}}{{\mathrm{d}}x^{k-1-r}}\,\bigl\{q(x)
S_{[k]}(x,Q)\bigr\} \delta q^{(r)}(x)|_-^+.\nonumber
\end{eqnarray}

At this point we restrict consideration to functions $\delta q$ whose
derivatives vanish sufficiently quickly at the end-points that we can
suppose the boundary terms in the last line of (\ref{eqfullvariationq})
vanish. Then (\ref{eqvarn}) will vanish for all such~$\delta q(\cdot)$
if
%
\begin{equation}\label{eqfirstvariationq}
\sum_{k=0}^{m}(-1)^{k+1}\,\frac{{\mathrm{d}}^{k}}{{\mathrm{d}}x^{k}}\,
\bigl\{q(x)
S_{[k]}(x,Q)\bigr\} \equiv\lambda_Q,
\end{equation}
that is, the left-hand side of (\ref{eqfirstvariationq}) is a~constant,
independent of $x$.

Motivated by (\ref{eqfirstvariationq}), we introduce the following
linear differential operator~$L$ on $q$-functions:
%
\begin{equation}
\label{eqlhat}
L \defeq\sum_{k\geq0}(-1)^{k+1}\D^k  q_0 \,\frac{\partial
}{\partial q_k}.
\end{equation}
Unless overridden by parentheses, operators here and elsewhere
associate to the right, so that $Tq_0 f$ means $T(q_0 f)$, that is, we have
%
\begin{equation}
\label{eqlhat0}
Lf = \sum_{k\geq0}(-1)^{k+1} \D^k\biggl(q_0  \,\frac{\partial
f}{\partial q_k}\biggr).
\end{equation}
For $f$ of order $m$, the series in (\ref{eqlhat0}) terminates at $k=m$,
and the order of~$L f$ is at most $2m$.

We can now write (\ref{eqfirstvariationq}) as
%
\begin{equation}
\label{eqlhatq}
L s = \lambda_Q,
\end{equation}
where equality in (\ref{eqlhatq}) is required to hold for all
$(x,q_0,q_1,\ldots,q_{2m})$ such that $q_j = q^{(j)}(x)$
($j=0,\ldots,2m$). In particular, a~\textit{sufficient} condition that
$S$ be weakly proper is that for some $\lambda\in\R$ we have
%
\begin{equation}
\label{eqlhatl}
L s = \lambda
\end{equation}
for all $x\in{\cal X}$, $\bldq \in\bQ_{2m}$.

So long as $\mathcal{P}$ is sufficiently large, the form (\ref{eqlhatl})
will also be necessary for (\ref{eqlhatq}) to hold. In particular,
suppose we impose the following condition on $\cal P$:

\begin{cond}
\label{cond2shot}
Given distinct $x_1, x_2 \in{\cal X}$, and any $\bldq
_1, \bldq _2 \in
\bQ_{2m}$, there exists $Q\in{\cal P}$ satisfying $q^{(j)}(x_1) =
q_{1,j}$, $q^{(j)}(x_2) = q_{2,j}$ $(j=0, \ldots,2m)$.
\end{cond}

Take arbitrary $Q_1, Q_2\in{\cal P}$, $x_1 \neq x_2\in{\cal X}$, and
set $q_{i,j} \defeq q_i^{(j)}(x_i)$ ($i=1,2;\allowbreak j=0, \ldots,2m$). Let
$Q$ be as given by Condition~\ref{cond2shot}. Evaluating (\ref
{eqlhatq}) at
$(x_1, \bldq _1)$ yields $\lambda_{Q_1} = \lambda_Q$,
and similarly
$\lambda_{Q_2} = \lambda_Q$. Thus $\lambda_Q$ cannot depend on $Q$,
and so (\ref{eqlhatl}) must hold. Moreover, taking $x_1 = x$, $
\bldq_1 =
\bldq $, and $x_2$, $\bldq _2$ arbitrary,
this must hold for any $x\in{\cal
X}$, $\bldq \in\bQ_{2m}$.

So we henceforth restrict attention to solutions of
(\ref{eqlhatl}). We note that a~particular solution of
(\ref{eqlhatl}) is given by the \textit{log-score}:
\[
s = -\lambda \ln q_0.
\]
Since $L$ is a~linear operator, the general solution is of the form
\[
s = -\lambda \ln q_0 + s_0,
\]
where $s_0$ satisfies the \textit{key equation}:
%
\begin{equation}
\label{eqlhatl0}
{L} s_0 = 0.
\end{equation}
Because of this we shall confine attention to solutions of the key
equation, and shall term any solution of (\ref{eqlhatl0}) a~\textit{key
local score function}.

\subsection{Connection to classical calculus of variations}
\label{secconnect}

Because the Lagrange multiplier $\lambda$ associated with a~key local
scoring rule $s$ vanishes, setting $q(x) \equiv p(x)$ will in fact
deliver a~\textit{globally} stationary point [i.e., without imposing the
normalization constraint $\int{\mathrm{d}} x\,q(x) = 1$] of the
corresponding expected score
\[
\int {\mathrm{d}}x\, p(x)
s\bigl\{x,q(x),q'(x),q''(x),\ldots,q^{(m)}(x)\bigr\}.
\]

The classical calculus of variations---see, for example,
van Brunt [(\citeyear{vanbrunt04}), equations (2.9), (3.3)]---would (again, ignoring the
boundary terms) identify the solution to this unconstrained
variational problem in $q(\cdot)$ as solving the \textit{Euler--Lagrange
equation}:
%
\begin{equation}
\label{eqel}
\Lambda  p_0 s(x, q_0, \ldots, q_m) = 0,
\end{equation}
where $\Lambda$ is the \textit{Lagrange operator}:
%
\begin{equation}
\label{eqLprime}
\Lambda\defeq\sum_{k\geq0}(-1)^{k}\D^k\,\frac{\partial}{\partial q_k}.
\end{equation}
We want the solution of (\ref{eqel}) to be $\bldq =
\bp$.

Now when evaluated at $\bldq = \bp$, $\Lambda  q_0 s
= s + \Lambda
p_0 s$. So $\bldq = \bp$ should satisfy
%
\begin{equation}
\label{eqneweq}
(I - \Lambda q_0)s = 0,
\end{equation}
where $I$ is the identity operator.\vadjust{\goodbreak}

But
\[
\frac{\partial}{\partial q_0} (q_j f)=
\cases{
q_j \,\dfrac{\partial f}{\partial q_0}, &\quad $(j>0)$,\vspace*{2pt}\cr
f +q_0\,\dfrac{\partial f}{\partial q_0}, &\quad $(j=0)$,}
\]
so we have
%
\begin{eqnarray}
\label{eqldashl}
L &=& I +\sum_{k\geq0}(-1)^{k+1}\D^k\,\frac{\partial}{\partial
q_k}q_0{\circ}\nonumber\\[-8pt]\\[-8pt]
&=& I - \Lambda q_0{\circ} .\nonumber
\end{eqnarray}
(Here and throughout, for $g$ a~$q$-function, $g\circ$ denotes the
multiplication operator $f \mapsto g f$, the optional symbol $\circ$
being attached, where required, to avoid confusion with the
$q$-function $g$ itself.)

Hence (\ref{eqneweq}) becomes $L s=0$, so recovering the key equation
(\ref{eqlhatl0}).

\section{Properties of differential operators}
\label{secsolve}
For a~further study of the key equation and its properties we shall
have a~detailed look at the differential operators introduced earlier,
together with some new ones. We recall
\begin{eqnarray*}
\D&=& \frac{\partial}{\partial x} + \sum_{j\geq0} q_{j+1}
\,\frac{\partial}{\partial q_j},\\
{L} &=& \sum_{k\geq0}(-1)^{k+1}\D^k  q_0 \,\frac{\partial}{\partial
q_k},\\
\Lambda&=& \sum_{k\geq0}(-1)^{k}\D^k\,\frac{\partial}{\partial
q_k}\\
&=& (I - L) q_0^{-1} {\circ}.
\end{eqnarray*}

\begin{lem}
\label{lemBoperatorD}
The Lagrange operator $\Lambda$ annihilates the total derivative
operator $\D$:
%
\begin{equation}
\label{eqldashd}
\Lambda\D= 0.
\end{equation}
\end{lem}
\begin{pf}
Using $(\partial/\partial q_k)D = D(\partial/\partial q_k) +
(\partial/\partial q_{k-1})$ we have
\begin{eqnarray*}
\Lambda D &=& \sum_{k\geq0} (-1)^{k} D ^{k+1}\,
\frac{\partial}{\partial q_k} + \sum_{k\geq0} (-1)^{k} D ^{k}\,
\frac{\partial}{\partial q_{k-1}}\\
&=& 0.
\end{eqnarray*}
\upqed\end{pf}

We now introduce the \textit{Euler operator}:
%
\begin{equation}
\label{eqm}
E \defeq\sum_{j\geq0} q_j \,\frac{\partial}{\partial q_j}.\vadjust{\goodbreak}
\end{equation}

\begin{lem}
\label{lemcommute}
The Euler operator $E$ commutes with $\D$ and with $L$, while
%
\begin{equation}
\label{eqlambdacommute}
\Lambda E = E\Lambda+ \Lambda.
\end{equation}
\end{lem}
\begin{pf}
From the easily verified relations
%
\begin{eqnarray}
\label{eqeasy1}
Eq_j\circ&=& q_j {\circ}+ q_j E,\nonumber\\[-8pt]\\[-8pt]
(\partial/\partial q_k)E &=& E(\partial/\partial q_k) +
(\partial/\partial q_k),\nonumber
\end{eqnarray}
it readily follows that $E$ commutes with $q_j(\partial/\partial
q_k)$. Since clearly $E$ commutes with $\partial/\partial x$, $E$
thus commutes with~$\D$, and consequently with any power of~$\D$.
From (\ref{eqlhat}), we now see that $E$ commutes with $L$.

Now (\ref{eqldashl}) gives that $E$ commutes with $\Lambda q_0{\circ}$,
and thus, noting from (\ref{eqeasy1}) that $q_0 E q_0^{-1} \circ= E -
I$, we have
\[
E \Lambda=
E\Lambda q_0q_0^{-1}\circ=
\Lambda q_0Eq_0^{-1}\circ=\Lambda(E - I)
=\Lambda E-\Lambda,
\]
which yields (\ref{eqlambdacommute}).
\end{pf}
%
\begin{theorem}
\label{thmml}
We have that $\Lambda E =\Lambda q_0 \Lambda$.
\end{theorem}
\begin{pf}
Using $q_k \D= \D q_k {\circ}- q_{k+1}\circ$, we can readily show by
induction that, for $k \geq0$,
\[
q_0 \D^k = \D\Biggl\{\sum_{j=0}^{k-1}(-1)^j q_j \D^{k-1-j}\Biggr\} +
(-1)^k q_k {\circ}.
\]
It now follows from (\ref{eqldashd}) that $\Lambda q_0 D^k = (-1)^k
\Lambda q_k \circ$. Applying this term-by-term to (\ref{eqLprime}) we
obtain $\Lambda q_0 \Lambda= \sum_k \Lambda q_k (\partial/\partial
q_k) = \Lambda E$.
\end{pf}

For later purposes we introduce, for any integer $r$,
%
\begin{equation}
\label{eqBoperator}
B_r \defeq\sum_{k\geq r+1}(-1)^{k-1-r}D^{k-1-r}\,
\frac{\partial}{\partial q_k}
\end{equation}
with the understanding ${\partial}/{\partial q_k} = 0$ if $k \leq0$
(in particular, $B_{-1} = \Lambda$). We further define
%
\begin{equation}
\label{eqBent}
C \defeq\sum_{r\geq0}  q_r B_r.
\end{equation}

\begin{lem}
\label{lemBoperatorDr}
It holds that
%
\begin{eqnarray}
\label{eqdb}
DB_r &=& (\partial/\partial q_r) - B_{r-1},\\
\label{eqbd}
B_r D &=& (\partial/\partial q_r).
\end{eqnarray}
\end{lem}
\begin{pf}
Equation (\ref{eqdb}) follows easily from the definition
(\ref{eqBoperator}), while (\ref{eqbd}) can be proved in the same
way as
Lemma~\ref{lemBoperatorD}.
\end{pf}
%
\begin{theorem} We have
\label{thmmcor}
%
\begin{eqnarray}
\label{eqcd}
CD &=& E,\\
\label{eqdc}
DC &=& E - q_0 \Lambda.
\end{eqnarray}
\end{theorem}
\begin{pf}
Equation (\ref{eqcd}) follows directly from (\ref{eqBent}) and
(\ref{eqbd}). From (\ref{eqdb}), $D q_r B_r = q_{r+1} B_r - q_r B_{r-1}
+ q_r (\partial/\partial q_r)$, and thus $DC = \sum_{r\geq0} q_r
(\partial/\partial q_r) - q_0 B_{-1} = E - q_0 \Lambda$.
\end{pf}

\section{Homogeneous scoring rules}
\label{sechomscore}

We shall see that all key local score functions, that is, solutions to the
key equation $Ls=0$, are homogeneous in the sense that changing
$\bldq$
by a~multiplicative factor does not change the value of $s$; hence the
associated scoring rule $S(x,Q)$ only involves the density $q$ up to a~constant factor. We shall formalize and show this below.
%
\begin{definer}
\label{defhom}
A~$q$-function $f$ is said to be \textit{homogeneous of degree $h$},
or \textit{$h$-homogeneous}, if, for any $\lambda> 0$, $f(x,
\lambda\bldq ) \equiv\lambda^h f(x,\bldq )$.
\end{definer}

With $E$ defined by (\ref{eqm}), Euler's homogeneous function theorem
implies that a~$q$-function $f$ is $h$-homogeneous if and only
%
\begin{equation}
\label{eqphistand}
Ef = h f.
\end{equation}

The partial derivatives $f_{[j]}={\partial f}/{\partial q_j}$ of an
$h$-homogeneous function are homogeneous of order $h-1$, while
$f_{[x]}={\partial f}/{\partial x}$ is homogeneous of order $h$. It
follows that, if $f$ is $h$-homogeneous, then so is $\D f$.

In this work we shall only need to deal with homogeneity of degree~$0$,
where $f(x, \lambda\bldq ) \equiv f(x,
\bldq)$, and of degree $1$,
where $f(x, \lambda\bldq ) \equiv\lambda f(x,
\bldq)$. Clearly, $f$ is
$0$-homogeneous if and only if $q_0 f$ is $1$-homogeneous.

A~scoring rule $S$ will be called \textit{homogeneous} if its score
function $s$ is $0$-homogeneous. As already noted, the logarithmic
score is $0$-local, but is not homogeneous. The Hyv\"arinen scoring
rule and its generalizations, as described in Section~\ref{sechyvscore}, are
$2$-local, and are homogeneous.

We can now easily show that key local score functions are homogeneous:
%
\begin{theorem}
\label{thmkey0}
If $Lf = 0$, then $f$ is $0$-homogeneous.
\end{theorem}
\begin{pf}
In this case $f=(I-L)f=\Lambda q_0f$ and thus from
(\ref{eqlambdacommute}) and Theorem~\ref{thmml} we get
\[
Ef = E\Lambda q_0 f= \Lambda E q_0 f-\Lambda q_0 f= \Lambda q_0
(\Lambda q_0 f) -\Lambda q_0 f = \Lambda q_0 f- \Lambda q_0 f=0
\]
as required.\vadjust{\goodbreak}
\end{pf}

We can further show that, if we consider the restriction of the
operator $L$ to $0$-homogeneous functions, it acts as a~\textit{projection
operator}; $I-L$ is then the complementary projection. This is a~consequence of the fact that $L$ is idempotent when restricted to
$0$-homogeneous functions:
%
\begin{theorem}
\label{thmkeyhom}
If $f$ is $0$-homogeneous, then so is $Lf$, and $L^2f = Lf$.
\end{theorem}
\begin{pf}
Since $E$ commutes with $L$, if $f$ is $0$-homogeneous we have $ELf
= LEf = 0$, so $Lf$ is $0$-homogeneous as well. If $f$ is
$0$-homogeneous, $q_0f$ is $1$-homogeneous and thus by
(\ref{eqphistand}) and Theorem~\ref{thmml}
\[
\Lambda q_0 f = \Lambda E q_0 f = \Lambda q_0 \Lambda q_0 f,
\]
so $I-L=\Lambda q_0 \circ$ is idempotent when restricted to
$0$-homogeneous functions, whence so is $L$.
\end{pf}

Elaborating the consequences of these results we get:
%
\begin{cor}
\label{corL0}
We have that:
%
\begin{longlist}
\item$Lf = 0$ if and only if $f= (I-L)g$ for some $0$-homogeneous
$g$; equivalently, $f = \Lambda\phi$ for some $1$-homogeneous
$\phi$;
\item
if $f$ is $0$-homogeneous, then $(I-L)f = 0$ if
and only if $f= Lg$ for some $0$-homogeneous $g$.
\end{longlist}
\end{cor}
\begin{pf}
If $Lf = 0$, then $f$ is $0$-homogeneous by Theorem~\ref{thmkey0}.
The other
properties are easy consequences of the fact that $L$ and $I-L$ are
complementary projections in the space of $0$-homogeneous functions.
\end{pf}

Collecting everything, we have the following main result:
%
\begin{theorem}\label{thmmain}
A~$q$-function $s$ is a~key local score function if and
only if any one (and then all) of the following conditions holds:
\begin{longlist}
\item The function $s$ satisfies the key equation $Ls=0$, where the
operator~$L$ is given by (\ref{eqlhat}).
\item
We can express
$s = (I-L)g$
where $g$ is a~$0$-homogeneous $q$-function.
\item
We can express
$s = \Lambda\phi$
where $\phi$ is a~$1$-homogeneous $q$-function and the operator
$\Lambda$ is given by (\ref{eqLprime}).
\end{longlist}
Moreover, $s$ is then $0$-homogeneous.
\end{theorem}

When (ii) above holds, we say that $s$ is \textit{derived from
$g$}; when (iii) holds, we say that $s$ is \textit{generated
by $\phi$}. The key local score function generated by a~$1$-homogeneous $q$-function $\phi$ of order $t$ is thus
\[
s(x, \bldq ) = \sum_{k=0}^t (-1)^{k} D^k \phi_{[k]}(x,
\bldq ).
\]
The only term in $s$ that involves $q_{2t}$ is $(-1)^{t}\phi_{[tt]}
q_{2t}$. In particular, if $\phi_{[tt]} \neq0$, $s$ is of exact
order $2t$. Hence we have demonstrated the existence of key local
scoring rules of all positive even orders.

The key local scoring rule $S$ generated by $\phi$ is then
%
\begin{equation}
\label{eqSmphi}
S(x, Q) = \sum_{k=0}^t (-1)^{k} \,\frac{{\mathrm{d}}^k }{{\mathrm{d}}
x^k}\,
\phi_{[k]}\bigl\{x, q(x), q'(x), \ldots, q^{(t)}(x)\bigr\}.
\end{equation}

For the case $t=1$ we obtain a~second-order rule:
\[
S(x, Q) = \phi_{[0]}\{x, q(x), q'(x)\} - \frac{{\mathrm{d}}} {{\mathrm{d}}
x}\, \phi_{[1]}\{x, q(x), q'(x)\},
\]
where $\phi(x, q_0, q_1)$ is $1$-homogeneous.

The Hyv\"arinen scoring rule (\ref{eqhyv1}) is generated in this way by
$\phi= -\half q_1^2/q_0$. More generally, choosing $\phi=
-q_1^k/q_0^{k-1}$ ($k\geq1$) yields
%
\begin{equation}
\label{eqo2}
S(x,Q) = (k-1)(y_1^k + k y_1^{k-2}y_2),
\end{equation}
where $y_i \defeq({\mathrm{d}}^i/{\mathrm{d}}x^i)\ln q(x)$. We can express a~general
$1$-homogeneous $x$-independent $q$-function of order $1$ as a~power
series:
%
\begin{equation}
\label{eqpower}
\phi(q_0,q_1) = q_0 \sum_{k\geq1} a_k (q_1/q_0)^k.
\end{equation}
Now combining the rules (\ref{eqo2}) arising from the individual terms
in (\ref{eqpower}), we obtain the series form of a~general
$x$-independent second-order scoring rule described by
\citet{ehmgneitingtr}.

\section{Gauge transformation}
\label{secnonunique}

The map $\phi\mapsto s = \Lambda\phi$ in Theorem~\ref{thmmain}(iii)
is many-to-one: two $1$-homogeneous functions $\phi_1$ and $\phi_2$
will generate the identical score function $s = \Lambda\phi_1 =
\Lambda\phi_2$ if and only if $\Lambda(\phi_2 - \phi_1) = 0$. And
this will hold if and only if $\phi_2 - \phi_1$ has the total
derivative form $D\psi$:
%
\begin{lem}
\label{lemconverse}
Suppose $\phi$ is $1$-homogeneous. Then $\Lambda\phi= 0$ if and
only if $\phi$ has the form $D\psi$.
\end{lem}
\begin{pf}
If $\phi= D\psi$, then $\Lambda\phi= 0$ by (\ref{eqldashd}).
Conversely, suppose $\phi$ is $1$-homogeneous and $\Lambda\phi=
0$. Then $q_0^{-1}\phi$ is $0$-homogeneous and $(I-L)q_0^{-1}\phi=
0$, so by Corollary~\ref{corL0}(ii) there exists
$0$-homogeneous $g$ such
that $ q_0^{-1}\phi= L g$. Now take $\psi= C q_0 g$, with $C$
given by (\ref{eqBent}). Then, using (\ref{eqdc}), $D\psi=
(E-q_0\Lambda)q_0 g = (I - q_0\Lambda)q_0 g$, since $Eq_0 g = q_0 g$
because $q_0 g$ is $1$-homogeneous; and this is $q_0(I - \Lambda
q_0\circ)g = q_0 L g = \phi$.
\end{pf}

Borrowing terminology from physics, we term a~transformation of the
form $\phi\rightarrow\phi+ D\psi$ a~\textit{gauge transformation};
the invariance of $s$ under such a~transformation of~$\phi$ is
\textit{gauge invariance}. The choice of a~particular function~$\phi$,
out of the equivalence class of functions differing only by a~total
derivative $D\psi$ and thus generating the same scoring rule, is a~\textit{gauge choice}.

Clearly if $\phi_2 - \phi_1 = D\psi$ and both $\phi_1$ and $\phi_2$
are $1$-homogeneous, then $D\psi$ must be $1$-homogeneous. This will
be so if $\psi$ is itself $1$-homogeneous. The converse also
essentially holds:
%
\begin{lem}
\label{lemdcon}
Suppose $D\psi$ is $1$-homogeneous. Then, for some constant~$a$,
$\psi+a$ is $1$-homogeneous.
\end{lem}
\begin{pf}
We have $ED\psi= D\psi$. Since by Lemma~\ref{lemcommute} $D$ commutes
with $E$, $D(E\psi- \psi) = 0$. Thus by Lemma~\ref{lemdint} $E\psi-
\psi$ is a~constant, $a$ say. Then $E(\psi+ a) = E\psi= \psi+a$,
so $\psi+ a$ is $1$-homogeneous.
\end{pf}

Since the addition of a~constant has no consequences for the analysis,
we henceforth call a~transformation $\phi\rightarrow\phi+ \kappa$ a~gauge transformation if and only if $\kappa$ has the form $D\psi$ with
$\psi$ $1$-homogeneous.

\subsection{Standard gauge choice}
\label{secstandard}
For any key local score function $s$ we note that
%
\begin{equation}
\label{eqstand}
\phi= q_0 s
\end{equation}
satisfies
\[
\Lambda\phi= \Lambda q_0 s=(I-L)s=s
\]
and hence $\phi= q_0 s$ is a~valid gauge choice for $s$. We call
(\ref{eqstand}) the \textit{standard gauge choice}.

\subsection{Equivalence}
\label{secequiv}
Suppose $s$ is generated by $\phi$, and let $\phi^* = \phi+ \chi$
with $\chi= a(x) q_0$. This is not a~gauge transformation if $a~\not
\equiv0$, but the score function it generates, $s^* = s + a(x)$, is
equivalent to $s$---which we describe by saying $\phi^*$ and $\phi$
are \textit{equivalent}. Conversely, if $\phi^\dag$ generates $s +
a(x)$ it must be a~gauge transformation of $\phi^*$, and hence of the
form $\phi+ a(x) q_0 + D\psi$---this form thus being necessary and
sufficient for equivalence. We note in particular that $\phi^\dag$ of
the form $\phi+ \sum_{k\geq0} a_k(x) q_k$ is equivalent to $\phi$,
since it generates $s^\dag= s + \sum_{k\geq0} (-1)^k a_k^{(k)}(x)$.

\subsection{Nonexistence of odd-order key local scores}
\label{secnonodd}

In Section~\ref{sechomscore} we established the existence of key local score
functions of all positive even orders. Here we show that no key local
score function can be of odd order.

Take $ s = \Lambda\phi$ as in Theorem~\ref{thmmain}(iii),
and suppose
$s$ has odd order. If $\phi$ is of order~$t$, the order of $s$ is at
most $2t$; since it is odd, it must be strictly less than~$2t$.
Again, the only term in $s$ that could possibly involve $q_{2t}$ is
$(-1)^{t}\phi_{[tt]} q_{2t}$, whence $\phi_{[tt]} = 0$. Hence
$\phi_{[t]}$, which must be $0$-homogeneous, has the form
$A(x$, $q_{0},\ldots,q_{t-1})$.\vadjust{\goodbreak}

Now define
%
\begin{equation}
\psi(x,q_{0},\ldots,q_{t-1}) \defeq
-\int_{0}^{q_{t-1}}A(x,q_{0},\ldots,q_{t-2},z) \,{\mathrm{d}}z
\end{equation}
[for the case $t=1$ the integrand on the right-hand side is $A(x,z)$].
It is easy to see that $\psi$ is $1$-homogeneous, and
%
\begin{equation}
\label{eqreduce}
\phi_{[t]} + \psi_{[t-1]} = 0.
\end{equation}

Let $\phi^* = \phi+ D\psi$, which is of order at most $t$. Since
this is a~gauge transformation, $\phi^*$ generates the same scoring
rule $s$ as $\phi$ does. But from (\ref{eqreduce}), $\phi^* _{[t]} =
\phi_{[t]} + \psi_{[t-1]} = 0$, so that $\phi^*$ is in fact of order
at most $t-1$, whence $s$ is of order at most $2t-2$. We can now
repeat the argument, stepping down~$t$ by 1 each time, until we reach
a~contradiction.

\subsection{Second-order rule}
\label{sec2ord}

A~similar argument to the above shows that, for any key local scoring
rule of exact even order $2t$, there exists a~gauge choice of exact
order $t$.

A~second-order rule can thus always be generated by a~$1$-homogeneous~$\phi$ of order $1$. However, a~change of gauge may increase the
order of the generating function---for example, the standard gauge
choice has order $2$.

If $\phi_1$ and $\phi_2$ are both gauge choices of order $1$, then
their difference is of order~$1$ and has the form $D\psi$ for some
$1$-homogeneous $\psi$. Then $\psi$ must be of order~$0$, and hence of
the form $\psi= c(x) q_0$. It follows that an order-$1$ gauge choice
is determined up to an additive term of the form $c'(x) q_0 +
c(x)q_1$. More generally, by Section~\ref{secequiv} two $1$-homogeneous
functions $\phi_1$ and $\phi_2$ of order $1$ are equivalent if their
difference has the linear form $a_0(x) q_0 + a_1(x) q_1$; and this is
also necessary, since, again by Section~\ref{secequiv}, $\phi_2$ must then
have the form $\phi_1 + a(x) q_0 + c'(x) q_0 + c(x)q_1$.

\section{Decomposition}
\label{secdecomp}

The variational analysis has identified the form (\ref{eqSmphi}), where
$\phi$ is $1$-homogeneous, for a~key local scoring rule. We now
consider the properties of such a~rule in more detail.

Starting from (\ref{eqSmphi}), we compute the expected score,
\begin{eqnarray*}
S(P,Q) &=& \int_-^+ {\mathrm{d}}x\, p(x) S(x,Q)\\
&=& \sum_{k\geq0} (-1)^{k} \int_-^+ {\mathrm{d}}x \,p(x) \,\frac{{\mathrm{d}}^k}
{{\mathrm{d}} x^k}\,
\phi_{[k]}\bigl\{x, q(x), q'(x), \ldots, q^{(t)}(x)\bigr\}
\end{eqnarray*}
by evaluating the $k$th term in the sum using the integration by parts
formula~(\ref{eqfullparts}). Collecting terms, we obtain
%
\begin{equation}
\label{eqcanonicalsplit}
S(P,Q)=S_0(P,Q)+S_+(P,Q)+S_-(P,Q),
\end{equation}
where the \textit{integral expected score} $S_0$ is given by
%
\begin{equation}
\label{eqexpectedorderall}
S_0(P,Q)= \int_-^+ {\mathrm{d}}x \sum_k p_k \phi_{[k]}(
\bldq)
\end{equation}
and
%
\begin{equation}
\label{eqpbcphi}
S_{\pm}(P,Q) = \mp S_b(\bp,\bldq )|_\pm,
\end{equation}
where the \textit{boundary expected score} $S_b$ is given by
%
\begin{equation}
\label{eqsb}
S_b(\bp,\bldq ) \defeq\sum_{r\geq0}  p_r B_r \phi
(\bldq )
\end{equation}
with $B_r$ defined by (\ref{eqBoperator}). In these formulas the
dependence on $x$ has been suppressed from the notation for
simplicity, and we interpret $p_k \defeq p^{(k)}(x)$, $q_k \defeq
q^{(k)}(x)$.

Correspondingly, the entropy $H(Q) = S(Q,Q)$ can be decomposed:
%
\begin{equation}
\label{eqentropysplit}
H(Q)=H_0(Q)+H_+(Q)+H_-(Q)
\end{equation}
with \textit{integral entropy}
%
\begin{equation}\label{entropy}
H_0(Q) \defeq\int_-^+ {\mathrm{d}}x \sum_k q_k  \phi_{[k]}(
\bldq)=
\int_-^+ {\mathrm{d}}x \,\phi(\bldq ),
\end{equation}
where the last equality follows from Euler's theorem (\ref{eqphistand});
and $H_\pm(Q) = \mp H_b(\bldq )|_\pm$,
where the \textit{boundary entropy} $H_b(\bldq )$ satisfies
\begin{eqnarray*}
H_b(\bldq ) &=& S_b(\bldq ,
\bldq)\\
&=& C \phi(\bldq )
\end{eqnarray*}
with the operator $C$ defined by (\ref{eqBent}).

The divergence now becomes
%
\begin{equation}
\label{eqdFfull}
d(P,Q)=d_0(P,Q)+d_+(P,Q)+d_-(P,Q),
\end{equation}
where $d_0(P,Q)=S_0(P,Q)-H_0(P)$, etc. In particular, the boundary
terms arise from the \textit{boundary divergence}
%
\begin{equation}
\label{eqbdiv}
d_b(P,Q) = \sum_r p_r B_r \{\phi(\bldq ) - \phi(\bp
)\}
\end{equation}
[where the final term involves substituting $\bp$ for
$\bldq$ after
computing $B_r \phi(\bldq )$]; while, using (\ref
{eqphistand}), the \textit{integral divergence} can be written as
%
\begin{equation}
\label{eqdF}
d_0(P,Q) = \int_-^+ {\mathrm{d}}x\,
\biggl[
\biggl\{
\phi(\bldq ) + \sum_k (p_{k} - q_{k})
\phi_{[k]}(\bldq )
\biggr\}
- \phi(\bp)
\biggr].
\end{equation}

It is easily seen that both $d_0$ and $d_b$ are unchanged by an
equivalence transformation $\phi^* = \phi+ \sum_{k\geq0} a_k(x)
q_k$.

\subsection{Change of gauge}
\label{secchangeg}

Although a~key local scoring rule $S$ is unchanged by a~gauge
transformation, the decompositions (\ref{eqcanonicalsplit}),
(\ref{eqentropysplit}) and (\ref{eqdFfull}), and in particular the
expression (\ref{eqdF}) for $d_0$, typically do change, terms being
redistributed between their constituents. Indeed, if we replace the
generating~$\phi$ by an alternative gauge choice
%
\begin{equation}
\label{eqphistar}
\phi^* = \phi+ D\psi,
\end{equation}
applying (\ref{eqexpectedorderall}) yields
\[
S_0^*(P,Q) = S_0(P,Q) + J
\]
with
\[
J \defeq\int_-^+  {\mathrm{d}}x \sum_k p_k
\biggl\{\frac{\partial}{\partial q_k}D\psi(
\bldq)\biggr\}.
\]
Using $({\partial}/{\partial q_k})D = D({\partial}/{\partial q_k}) +
{\partial}/{\partial q_{k-1}}$, and the interpretation of $D$ as
$\mathrm{d}/\mathrm{d}x$, this reduces to
\begin{eqnarray*}
J &=& \int_-^+  {\mathrm{d}}x \,\frac{\mathrm{d}}{{\mathrm{d}}x}
\sum_k p_k \psi_{[k]}(\bldq )\\
&=& \widehat S_+ + \widehat S_-,
\end{eqnarray*}
where $\widehat S_+\defeq\widehat S(\bp,\bldq
)|_+$,
$\widehat S_-\defeq-\widehat S(\bp,
\bldq)|_-$, with
$\widehat S(\bp,\bldq ) \defeq\sum_k p_k \psi
_{[k]}(\bldq )$.

Similarly, from (\ref{eqpbcphi}) and (\ref{eqsb}) we find the boundary
expected score transforming as
%
\begin{eqnarray}
S^*_b(\bp,\bldq ) &=& S_b(\bp,\bldq )+
\sum_{k} p_k B_k D \psi\nonumber\\[-8pt]\\[-8pt]
&=& S_b(\bp,\bldq )+ \widehat S(\bp,\bldq )
\nonumber
\end{eqnarray}
on using (\ref{eqbd}). The changes to the boundary terms thus
compensate exactly (as they must) for the changes to the integral
term.

We now have
\[
H_0^*(P) = H_0(P) + \widehat H_+ + \widehat H_-
\]
with $\widehat H_\pm\defeq\pm\widehat H|_\pm$ and
$\widehat H(\bp) = \psi(\bp)$; this follows from (\ref{eqcd}) since
$\psi$ is $1$-homogeneous. Correspondingly the boundary entropy
transforms as $H_b^*(\bp) = H_b + \psi(\bp)$.

It is notable that there is always a~gauge choice for which the
boundary entropy vanishes. Specifically:
%
\begin{theorem}
\label{thmbdent}
Let $s$ be a~key local score function. Then for the standard gauge
choice $\phi= q_0 s$, the boundary entropy function $H_b$ is
identically $0$.
\end{theorem}
\begin{pf}
From (\ref{eqdc}), $DH_b = DC\phi= E \phi- q_0 \Lambda\phi$. Since
$\phi$ is $1$-homogeneous and $s = \Lambda\phi$, this becomes $\phi
- q_0s = 0$. So $0 = C D H_b = E H_b$ by (\ref{eqcd}). But $EH_b =
H_b$ since $H_b$ is $1$-homogeneous.\vadjust{\goodbreak}
\end{pf}

The effect on a~gauge transformation on the decomposition of the
divergence is
%
\begin{equation}
\label{eqdchange}
d_0^*(P,Q) = d_0(P,Q) + \widehat d_+ + \widehat d_-,
\end{equation}
where $\widehat d_\pm\defeq\pm\widehat d|_\pm$ with
%
\begin{eqnarray}
\label{eqdadd}
\widehat d(\bp,\bldq ) &=& \sum_k p_k \psi
_{[k]}(\bldq )-\psi(\bp)\nonumber\\[-8pt]\\[-8pt]
&=& \psi(\bldq ) + \sum_k (p_k-q_k) \psi_{[k]}(
\bldq) -\psi(\bp)\nonumber
\end{eqnarray}
and with a~compensating change to the boundary divergence $d_b$.

\section{Propriety}
\label{secprop}

In this section we investigate the propriety of a~key local scoring
rule $S$. The scoring rule $S$ will be proper if and only if $d(P,Q)
\geq0$ for all $P$, \mbox{$Q \in{\cal P}$}. Clearly it is sufficient to
require nonnegativity of each term in the right-hand side of the
decomposition (\ref{eqdFfull}), and we proceed on this basis. We
investigate $d_+$ and $d_-$ in Section~\ref{secboundary} below; here we
consider the integral term $d_0$.

We note the similarity between formula (\ref{eqdF}) and that for the
Bregman divergence (\ref{eqbregd}) (especially where that is extended,
as in Section~\ref{secextbreg}, to allow further dependence of $\phi$ on $x$).
Correspondingly, concavity of the defining function plays a~crucial
role here, too.
%
\begin{definer}
\label{defconvex}
We call a~$1$-homogeneous $q$-function $\phi(x, \bldq)$
\textit{concave} if, for every $x\in\mathcal{X}$, $\bldq _1,
\bldq _2 \in
\bQ$,
%
\begin{equation}
\label{eqconvex}
\phi(x, \bldq _1 + \bldq _2) \leq\phi(x,
\bldq _1) + \phi(x, \bldq _2)
\end{equation}
(this is readily seen to be equivalent to the usual definition of
concavity in $\bldq $, for each $x$); and \textit{strictly
concave} if
strict inequality in (\ref{eqconvex}) holds whenever the vectors
$\bldq _1$ and $\bldq _2$ are linearly
independent.
\end{definer}
%
\begin{theorem}
\label{thmconvex}
Suppose that the scoring rule $S$ is generated by a~concave
$1$-homogeneous $q$-function $\phi$. Then $d_0(P,Q)$, as given by
(\ref{eqdF}), is nonnegative. Further, if $\phi$ is strictly concave,
then $d_0(P,Q)= 0$ if and only if $Q=P$.
\end{theorem}
\begin{pf}
Concavity implies that the integrand of (\ref{eqdF}) is nonnegative
for each~$x$; under strict concavity it will be strictly positive
with positive probability when $Q \neq P$.
\end{pf}
%
\begin{cor}
\label{corvanprop}
Suppose the conditions of Theorem~\ref{thmconvex} apply, and the boundary
terms $d_+(P,Q)$ and $d_-(P,Q)$ in (\ref{eqdFfull}) vanish identically
for \mbox{$P,Q \in{\cal P}$}. Then the (local, homogeneous) scoring rule
(\ref{eqSmphi}) is proper (strictly proper if $\phi$ is strictly
concave).
\end{cor}

\subsection{Checking concavity}

Given a~$1$-homogeneous $q$-function $\phi$ of order~$m$, define, for
$\bu= (u_1,\ldots,u_m)\in\R^m$,
%
\begin{equation}
\label{eqphitoF}
\Phi(x, \bu) \defeq\phi(x, 1, \bu).
\end{equation}
Then $\phi(x, \bldq )$ is determined by $\Phi$:
%
\begin{equation}
\label{eqFtophi}
\phi(x,\bldq ) = q_0 \Phi(x,\bu)
\end{equation}
with $u_i = {q_i}/{q_0}$ ($i\geq1$). It is often easier to check
concavity for $\Phi$ than for~$\phi$, and this is enough:
%
\begin{lem}
\label{lemconvF}
$\Phi$ is concave in $\bu$ if and only if $\phi$ is concave in
$\bldq $.
\end{lem}
\begin{pf}
``If'' follows immediately from (\ref{eqphitoF}). Conversely, if
$\Phi$ is concave,
\begin{eqnarray*}
\phi(x, \bp+ \bldq ) &=&
( p_0 + q_0) \Phi\biggl(x, \frac{p_0}{p_0+q_0} \frac{\bp}{p_0} +
\frac{q_0}{p_0+q_0} \frac{\bldq }{q_0}\biggr)\\
&\geq& p_0 \Phi\biggl(x, \frac{\bp}{p_0}\biggr)
+ q_0 \Phi\biggl(x, \frac{\bldq }{q_0}\biggr)\\
&=& \phi(x, \bp) + \phi(x, \bldq ).
\end{eqnarray*}
\upqed\end{pf}

It is further easy to see that $\Phi$ is strictly concave in $\bu$, in
the usual sense, if and only if $\phi$ is strictly concave in $\bldq$
in the sense of Definition~\ref{defconvex}.

\subsection{Change of gauge}
\label{secchangeg2}

Even if the initial gauge choice $\phi$ is concave in $\bldq$, so that
$d_0(\bp,\bldq ) \geq0$, under a~gauge transformation (\ref{eqphistar})
the term $\widehat d(\bp,\bldq )$, as given by~(\ref {eqdadd}), means
that the gauge-transformed integral divergence term~$d^*_0$, given
by~(\ref{eqdchange}), need not be nonnegative; this would hold if the
resulting gauge choice $\phi^*$ were itself concave, but typically this
will not be so.

Note that if $\psi$ in (\ref{eqphistar}) is concave, then $\widehat
d(\bp,\bldq )\geq0$. However, this does not ensure
positivity of both
additional terms, since while the added term $\widehat d_+ =
\widehat d(\bp,\bldq )|_+$ will then be
nonnegative, the
other added term $\widehat d_- = - \widehat d(\bp,\bldq)|_-$
will be nonpositive.
%
\begin{expl}
\label{exgauge} The Hyv\"{a}rinen scoring rule (\ref{eqhyv1}) on
${\cal X} = \R$ is generated by the strictly concave $q$-function
$\phi=-\frac{1}{2}q_1^2/q_0$. Using this gauge choice in~(\ref{eqdF})
yields [cf. (\ref{eqhyvd})]
%
\begin{equation}
\label{eqhyv0}
d_0(P,Q) = \half\int  {\mathrm{d}}x\, p(x) (v_1-u_1)^2
\end{equation}
with $u_i \defeq q^{(i)}(x)/q(x)$, $v_i \defeq p^{(i)}(x)/p(x)$.

Alternatively we might use the standard gauge choice, $q_2 -
\frac{1}{2}q_1^2/q_0$, which is also strictly concave, and indeed
yields the same expression (\ref{eqhyv0}).

Now let $\psi\defeq- \half q_1 \ln(q_1/q_0)$, so that $D\psi=
-\half\{q_2 \ln(q_1/q_0) + q_2 - q_1^2/q_0\}$. Then $\phi^* = \phi
+ D\psi= -\half q_2 \{1 + \ln(q_1/q_0)\}$ is another possible gauge
choice, generating the identical scoring rule $S$. However,
$\phi^*$ is not concave,\vadjust{\goodbreak} and the integral divergence term (\ref{eqdF})
associated with $\phi^*$ is
\[
d^*_0(P,Q)=\half\int {\mathrm{d}}x\, p(x)\biggl\{u_2\biggl(1-
\frac{v_1}{u_1}\biggr) + v_2\ln\frac{v_1}{u_1}\biggr\},
\]
which is not nonnegative. In this case the extra terms in
(\ref{eqdchange}) arise from $\widehat d(\bp,\bldq ) =
\half p_0\{u_1 -
v_1 + v_1\ln(v_1/u_1)\}$.
\end{expl}

In the light of the above example it might be conjectured that, if $s$
can be generated from \textit{some} concave gauge choice, then the
standard gauge choice $\phi= q_0 s$ will be concave---equivalently,
from Lemma~\ref{lemconvF}, $s$ itself will be a~concave function of
the $u_i
= {q_i}/{q_0}$ ($i\geq1$)---but this need not hold:
%
\begin{expl}
\label{exstandgauge}
Take $\Phi= -u_1^4$ in (\ref{eqFtophi}). Then $\Phi$, and hence
$\phi$, is concave, but $s = 12 u_1^2 u_2 - 9 u_1^4$ is not concave.
\end{expl}

\section{Boundary issues}
\label{secboundary}

The boundary divergence terms in (\ref{eqdFfull}) are $d_\pm(P,\allowbreak Q) =
\mp
d_b(\bp,\bldq )|_\pm$, where $d_b$ is given by (\ref
{eqbdiv}). Their
behavior will depend on the family ${\cal P}$ of distributions under
consideration, and specifically on the behavior, at the end-points $+$
and $-$, of the densities of distributions in ${\cal P}$.

For propriety of these terms, we want $d_b(\bp,\bldq )$
to be positive at
the lower end-point $-$, and negative at the upper end-point $+$, for
all $P, Q \in{\cal P}$. For simplicity we might impose conditions on
${\cal P}$ sufficient to ensure that, for all densities $p(\cdot)$,
\mbox{$q(\cdot) \in{\cal P}$}, $d_b(\bp,\bldq )$ vanishes at
the end-points. A~family ${\cal P}$ having this property may be termed \textit{valid}
(with respect to the generating function $\phi$). However, there does
not appear to be a~natural choice for such a~valid class ${\cal P}$.
In particular, if ${\cal P}$ and ${\cal P}'$ are both valid families,
it does not follow that their union will be.

Note that the validity requirement depends on the gauge choice $\phi$,
and a~change of gauge could assist in ensuring that it holds.

For the special case of the standard gauge choice, $\phi^* = q_0 s$,
we know from Theorem~\ref{thmbdent} that the boundary entropy $H^*_b$
vanishes. \textit{If} the boundary quantities $S^*_b$, $H^*_b$,
$d^*_b$ behaved like regular quantities $S$, $H$, $d$ we could deduce
$d^*_b = 0$ [\citet{apd94}]; but this is a~big ``if,'' and the result
will not hold without imposing further conditions.

\subsection{Second-order rules}
\label{sec2prop}
For a~second-order rule with $1$-{homogeneous} generator $\phi(x, q_0,
q_1)$, we find
\[
d_b = p_0 \bigl\{\phi_{[1]}(\bldq ) - \phi_{[1]}(\bp)\bigr\}.
\]
Alternatively, the standard gauge choice is $\phi^* = \phi+ D\psi$
with $\psi= -C\phi= -q_0\phi_{[1]}$. From Section~\ref{secchangeg}, we find
\[
d^*_b(\bp,\bldq ) = S^*_b(\bp,\bldq ) =
-q_0 \bigl(p_0 \phi_{[01]} + p_1
\phi_{[11]}\bigr).
\]
That this vanishes (as we know from Theorem~\ref{thmbdent} it must)
for $\bp=
\bldq $ may be seen\vadjust{\goodbreak} on differentiating the relation
$\phi= q_0
\phi_{[0]} + q_1 \phi_{[1]}$ with respect to $q_1$; that it does not
depend on the choice of gauge $\phi$ of order $1$ follows from
Section~\ref{sec2ord}.

With $p_i = p^{(i)}(x)$, etc., we want $d_b(\bp, \bldq)$ [or, for the
standard gauge choice, $d^*_b(\bp, \bldq )$] to vanish
in the limit as we
approach the end-points $-$ and $+$, for all densities $p(\cdot)$,
$q(\cdot)$ of distributions in ${\cal P}$. Conditions for validity
will thus involve the behavior of $p(x)$ and $p'(x)$ at these
end-points.

For example, for the Hyv\"arinen rule, with gauge choice $\phi=
-\half q_1^2/q_0$, we require
%
\begin{equation}
\label{eqhyvdb}
d_b(\bp,\bldq ) = p_0\biggl(\frac{p_1}{p_0} - \frac
{q_1}{q_0}\biggr)
\rightarrow0
\end{equation}
as we approach the end-points of ${\cal X}$. (The same expression for
$d_b$ arises if we use the standard gauge choice $\phi^* = q_2 -\half
q_1^2/q_0$, which in this case is equivalent to~$\phi$.) To ensure
(\ref{eqhyvdb}) we might require, for example, that, for all densities
$p(\cdot)$ in~${\cal P}$, $\lim_{x \rightarrow\pm} p(x) = 0$ and
$\lim_{x \rightarrow\pm} p'(x)/p(x)$ is finite. However, this
excludes the possibility that both $p$ and $q$ are normal densities on
${\cal X} = \R$, even though, with this choice, $d_b$ as given by
(\ref{eqhyvdb}) does vanish at $\pm\infty$.
Ehm and Gneiting (\citeyear{ehmgneitingtr,ehmgneitingorder2}) described alternative
conditions on ${\cal P}$ that do admit this case.

In the ideal situation we will have a~(strictly) concave
$1$-homogeneous $\phi$, and a~family ${\cal P}$ valid with respect to
$\phi$. Then the associated key local scoring rule $S$ will be
(strictly) proper.

\section{Transformation of the data}
\label{sectransform}

So far we have considered a~variable $X$ taking values in a~real
interval ${\cal X}$, and have made essential use of the Euclidean
structure of ${\cal X}$ to define probability densities, derivatives,
etc. Taking a~step backward, suppose we start with an abstract
measurable sample space (the \textit{basic sample space}) ${\cal X^*}$,
a~\textit{basic variable} $X^*$ taking values in ${\cal X}^*$, and a~collection ${\cal P}^*$ of \textit{basic distributions} for $X^*$ over
${\cal X}^*$. Without assuming any further structure, we can define a~\textit{basic scoring rule} $S^*\dvtx {\cal X}^* \times{\cal P}^*
\rightarrow\reals$, and introduce the property of (strict) propriety,
exactly as before. However, at this level of generality it is less
straightforward to define what we should mean by saying that a~basic
scoring rule is \textit{local}. To do this we proceed as follows.

We suppose given a~collection $\Xi= \{\xi\}$ of \textit{charts}, where
each $\xi$ is an invertible measurable function from ${\cal X}^*$ onto
some open interval ${\cal X} \subseteq\reals$, and such that, for
$\xi, \overline\xi\in\Xi$, the composition $\overline\xi\xi
^{-1}\dvtx{\cal
X}\rightarrow\overline{\cal X}$ is \textit{smooth} and \textit{regular},
that is, infinitely often differentiable with strictly positive first
derivative. In other words, the basic space is a~one-dimensional
simply connected smooth manifold.

Picking any specific chart $\xi$ produces a~concrete
\textit{representation} of the abstract basic structure, in terms of the
real variable $X \defeq\xi(X^*)$, and, for any $Q^*\in{\cal P}^*$,
the induced distribution $Q$ for $X$ on ${\cal X} \subseteq\reals$
[so that $Q(A) = Q^*\{\xi^{-1}(A)\}$]; we take ${\cal P} \defeq\{P\dvtx
P^* \in{\cal P}^*\}$. Correspondingly, a~\textit{basic function}
$f^*\dvtx{\cal X}^* \times{\cal P}^* \rightarrow\reals$ (e.g., a~scoring
rule) is \textit{represented} by $f\dvtx{\cal X} \times{\cal
P}\rightarrow\reals$, such that $f(x,Q) = f^*(x^*, Q^*)$.

Let $\xi, \overline\xi$ be two such charts, and $X = \xi(X^*)$,
$\overline X = \overline\xi(X^*)$, etc. Then $\overline X =
\gamma(X)$, where $\gamma= \overline\xi\xi^{-1}$ is strictly
increasing, and both $\gamma$ and $\delta\defeq\gamma^{-1}$ are
smooth and regular. A~given basic distribution $Q^*$ for $X^*$ can be
represented \textit{either} by the distribution $Q$, for $X$,
\textit{or} by $\overline Q$, for $\overline X$. We assume that $Q$ has a~density function, $q(\cdot)$, with respect to Lebesgue measure on
${\cal X}$; then the density function $\overline q(\cdot)$ of
$\overline Q$ with respect to Lebesgue measure on $\overline{\cal X}$
will likewise exist, and, with $\overline x = \gamma(x)$, we will have
%
\begin{equation}
\label{eqqbarq}
\overline q(\overline x) = q(x)\,\frac{\mathrm{d}x}{\mathrm{d}\overline x} =
\alpha(x)  q(x)
\end{equation}
with $\alpha(x): = \gamma'(x)^{-1}$. An easy induction shows that we
can express
%
\begin{equation}
\label{eqexplicit2}
\overline q^{(k)}(\overline x) = \overline T_k\bigl( x, q(x),
\ldots,q^{(k)}(x)\bigr),
\end{equation}
where $\overline T_k$ has the form
%
\begin{equation}
\label{eqtform}
\overline T_k(x, q_0, \ldots, q_k) = \sum_{r} a_{kr}(x)  q_r
\end{equation}
and the coefficients $a_{kr}(x)$ satisfy $a_{kr}(x) = 0$ unless $0
\leq r\leq k$, $a_{00}(x) = \alpha(x)$, and
%
\begin{equation}
\label{eqarecur}
a_{k+1,r}(x) = \alpha(x) \{a'_{kr}(x) + a_{k,r-1}(x)\}.
\end{equation}
In similar fashion we can express
%
\begin{eqnarray}
\label{eqexplicit1}
q^{(k)}(x) &=& T_k\bigl(\overline x,
\overline q(\overline x), \ldots, \overline q^{(k)}(\overline
x)\bigr)\nonumber\\[-8pt]\\[-8pt]
&=& \sum_{r} \overline a_{kr}(\overline x)
\overline q^{(r)}(\overline x).\nonumber
\end{eqnarray}

It readily follows from (\ref{eqexplicit2}) and (\ref{eqexplicit1})
that a~basic function $f^*(x^*,Q^*)$ can be written, in the
$\xi$-representation, in the form $f(x, q(x), q'(x),\ldots,
q^{(m)}(x))$ if and only if the analogous property holds in the
$\overline\xi$-representation: $f^*(x^*,Q^*) = \overline f(\overline
x, \overline q(\overline x), \overline q'(\overline x),\ldots,
\overline q^{(m)}(\overline x))$. That is, the property of being
\textit{$m$-local} is independent of the particular representation used.
When this property holds for one, and thus for all, representations,
we can say that the basic function $f^*(x^*, Q^*)$ itself is
\textit{$m$-local}; a~$q$-function $f$ such that $f^*(x^*,Q^*) = f(x,
q(x),\allowbreak q'(x),\ldots, q^{(m)}(x))$ is the \textit{$\xi$-representation}
of $f^*$. We denote the vector space of all local basic functions by
${\cal V}^*$.

At a~more abstract level, motivated by (\ref{eqexplicit2}) and
(\ref{eqtform}), we define variables
%
\begin{eqnarray}
\label{eqqbar}
\overline x &\defeq&
\gamma(x), \nonumber\\[-8pt]\\[-8pt]
\overline q_k &\defeq& \overline T_k(x, q_0, \ldots
,q_k)
= \sum_{r} a_{kr}(x)  q_r.\nonumber
\end{eqnarray}
Inversely, we will then have
%
\begin{eqnarray}
\label{eqdefq}
x &=& \delta(\overline x), \nonumber\\[-8pt]\\[-8pt]
q_k &=& T_k(\overline x, \overline q_0, \ldots, \overline q_k)
= \sum_{r} \overline a_{kr}(\overline x)  \overline q_r.\nonumber
\end{eqnarray}

Using (\ref{eqdefq}), any $q$-function of order $m$, $f(x,q_0, \ldots,
q_m)$ can be rewritten as $\overline f(\overline x, \overline q_0,
\ldots, \overline q_m)$. If $f^*\in{\cal V}^*$ has $\xi$- and
$\overline\xi$-representations $f$ and $\overline f$, respectively,
then $\overline f$ can be obtained by reexpressing $f$ in this way.
Since $T_k$ is $1$-homogeneous, $f$ is homogeneous of degree $h$ in
the $q$'s if and only if $\overline f$ is homogeneous of degree $h$ in
the $\overline q$'s. In this case we may term the underlying local
basic function $f^*\in{\cal V}^*$ \textit{$h$-homogeneous}. Likewise,
since (for fixed $x$ or~$\overline x$) the functions $T_k$ and
$\overline T_k$ are linear, $f$ is (strictly) concave in the $q$'s if
and only if $\overline f$ is (strictly) concave in the $\overline
q$'s---in which case we may term $f^*$ itself (\textit{strictly})
\textit{concave}.

\subsection{Invariant operators}
\label{secinvop}

The linear differential operators $D$ and $L$ have only been defined
in terms of a~specific representation of the problem on the real line,
as determined by some chart $\xi$. Applying these definitions
starting from a~different real representation, determined by a~chart
$\overline\xi$, we will obtain possibly different operators,
$\overline D$, $\overline L$. The following results relate these. We
need the following lemma:
%
\begin{lem}
\label{lempartials}
We have
%
\begin{eqnarray}
\label{eqdiffq}
\frac{\partial}{\partial q_r} &=& \sum_{k} a_{kr}\,
\frac{\partial}{\partial\overline q_k},\\
\label{eqdxbar}
\frac{\partial}{\partial x} &=&
\alpha^{-1}\,\frac{\partial}{\partial\overline x} + \sum_r q_r\sum_k
a'_{kr} \,\frac{\partial}{\partial\overline q_k}.
\end{eqnarray}
\end{lem}
\begin{pf}
Equation (\ref{eqdiffq}) follows immediately from (\ref{eqqbar}). For
(\ref{eqdxbar}) we have
\[
\frac{\partial}{\partial x} =
\frac{\mathrm{d}\overline x}{\mathrm{d}x} \,\frac{\partial}{\partial\overline x} +
\sum_k \frac{\partial\overline q_k}{\partial x}\,
\frac{\partial}{\partial\overline q_k}.
\]
But ${\mathrm{d}\overline x}/{\mathrm{d}x} = \alpha^{-1}$, while from (\ref{eqqbar})
\[
\frac{\partial\overline q_k}{\partial x} = \sum_r a'_{kr} q_r,
\]
so (\ref{eqdxbar}) follows.
\end{pf}

We now show that if $f$ and $\overline f$ are, respectively, the $\xi$
and $\overline\xi$ representations of the same basic function $f^*$,
then $\overline D  \overline f$ is the $\overline\xi$-representation
of the basic function whose $\xi$-representation is $\alpha(x) D f$.
Note that the function $\alpha$, and hence the basic function so
represented, will depend on the charts considered.

\begin{theorem}
\label{thmdtrans} It holds that
\[
\overline D = \alpha(x) D.
\]
\end{theorem}

\begin{pf}
Informally, we observe that $D$ corresponds to the total derivative
$\mathrm{d}/\mathrm{d}x$ and $\overline D$ to $\mathrm{d}/\mathrm{d}\overline x$. Thus we expect
$\overline D = (\mathrm{d}x/\mathrm{d}\overline x) D$.

More formally, we have
\begin{eqnarray*}
\overline D &=&
\frac{\partial}{\partial\overline x}
+ \sum_{k} \overline q_{k+1} \,\frac{\partial}{\partial\overline
q_k}\\
&=& \frac{\partial}{\partial\overline x}
+ \sum_r q_r\sum_k
a_{k+1,r} \,\frac{\partial}{\partial\overline q_k}
\end{eqnarray*}
on using (\ref{eqqbar}). From (\ref{eqarecur}) this is
\[
\frac{\partial}{\partial\overline x}
+ \alpha\sum_r q_r\sum_k
(a'_{k,r} + a_{k,r-1}) \,\frac{\partial}{\partial\overline q_k}.
\]
On applying Lemma~\ref{lempartials} this reduces to $\alpha D$.
\end{pf}

Since, by the transformation rule (\ref{eqqbarq}) for densities,
$\overline q_0 = \alpha(x) q_0$, we thus have
%
\begin{cor}
\label{cordtrans}It holds that $\overline q_0^{-1}\overline D =
q_0^{-1}D$.
\end{cor}

It follows from Corollary~\ref{cordtrans} that, for $f^*\in{\cal
V}^*$, there
exists $g^*\in{\cal V}^*$ such that, in any representation, $g =
q_0^{-1}Df$. This shows the existence of an ``invariant'' linear
operator $D^*$ on ${\cal V}^*$ such that, in any representation,
$D^*f^*$ is represented by $q_0^{-1}Df$.

We next show that there exists an invariant linear operator $L^*$ on
${\cal V}^*$ such that, in any representation, if $f^*$ is represented
by $f$, then $L^*f^*$ is represented by $Lf$.
%
\begin{theorem} We have
\label{thmltrans}
$\overline L = L$.
\end{theorem}
\begin{pf}
On substituting (\ref{eqqbarq}) and (\ref{eqdiffq}) into the definition
(\ref{eqlhat}) of $L$ and rearranging, we obtain
%
\begin{equation}
L = \sum_k (-1)^{k+1} A_k \alpha^{-1} \overline q_0\,
\frac{\partial}{\partial\overline q_k},
\end{equation}
where the operator $A_k$ is given by
\[
A_k = \sum_r (-1)^{k-r} D^r a_{kr}{\circ}.
\]
The theorem will thus be proved if we can show $A_k = \overline
D^k\alpha\circ$; that is, using Theorem~\ref{thmdtrans}, we have to show:
%
\begin{equation}
\label{eqarind}
H_k\dvtx (\alpha D)^k\alpha\circ= \sum_r (-1)^{k-r} D^r a_{kr}{\circ}.
\end{equation}
We prove (\ref{eqarind}) by induction on $k$. First, $H_0$ holds
since both sides reduce to $\alpha\circ$. Now suppose $H_k$ holds.
Then
%
\begin{eqnarray}
\label{eqaa1}
(\alpha D)^{k+1}\alpha\circ&=& (\alpha D)^{k}\alpha D
\alpha{\circ}\nonumber\\[-8pt]\\[-8pt]
&=& \sum_r (-1)^{k-r} D^r a_{kr} D \alpha{\circ}.\nonumber
\end{eqnarray}
But $ a_{kr}D = (D a_{kr}\circ) - a'_{kr}\circ$, so that (\ref{eqaa1})
becomes
\[
\sum_r (-1)^{k-r} D^{r+1} a_{kr}\alpha{\circ}
- \sum_r (-1)^{k-r} D^r a'_{kr} \alpha{\circ},
\]
which can be written as
\[
\sum_r (-1)^{k + 1-r} D^r \alpha(a'_{kr} + a_{k,r-1}
)\circ
\]
and on applying (\ref{eqarecur}) we have verified $H_{k+1}$.
\end{pf}

\subsection{Invariance of scoring rule}
\label{secinvarscore}

On applying Theorem~\ref{thmltrans}, we see that the general
homogeneous key
local scoring rule, as given by (ii) of Theorem~\ref{thmmain}, can be
expressed invariantly as
\[
S^*(x^*, Q^*) = (I-L^*)g^*(x^*,Q^*),
\]
where $g^*$ is a~$0$-homogeneous local basic function. Then, in any
representation, we will have $S(x,Q) = (I-L)g(x,Q)$. We may thus say
that the scoring rule $S^*$ is \textit{derived from} the local basic
function $g^*$. In particular the expected score $S^*(P^*,Q^*)$, and
consequently the entropy function $H^*(P^*)$ and the divergence
function $d^*(P^*,Q^*)$, are fully determined by the basic function
$g^*$, independently of how that may be represented.

In fact more is true: the individual components $S_0(P,Q)$,
$S_+(P,Q)$,\break $S_-(P,Q)$ of $S(P,Q)$, in the decomposition
(\ref{eqcanonicalsplit}) arising from the integration by parts, each
correspond to an invariant expression $S_0^*(P^*,Q^*)$,
$S_+^*(P^*,Q^*)$, $S_-^*(P^*,Q^*)$ (and similarly for the
decompositions of $H$ and $d$).

We show this first for the integral term $S_0$. We need the following
lemma, showing that the expression $\Pi\defeq\sum_r p_r
\, {\partial}/{\partial q_r}$ represents an invariant operator $\Pi^*$
(depending on a~distribution $P^*$, and acting on a~function of
a~distribution $Q^*$, both defined over ${\cal V}^*$).

\begin{lem}
\label{lempq} We have
\[
\sum_r p_r \,\frac{\partial}{\partial q_r}
= \sum_k \overline p_k \,\frac{\partial}{\partial\overline q_k}.\vadjust{\goodbreak}
\]
\end{lem}
\begin{pf}
Using (\ref{eqdiffq}) and (\ref{eqqbar}), we have
\begin{eqnarray*}
\sum_r p_r \,\frac{\partial}{\partial q_r} &=&
\sum_k \frac{\partial}{\partial\overline q_k}
\sum_{r=0}^k a_{kr}(x)   p_r\\
&=& \sum_k \overline p_k \,\frac{\partial}{\partial\overline q_k}.
\end{eqnarray*}
\upqed\end{pf}

Now consider expression (\ref{eqexpectedorderall}) for $S_0(P,Q)$,
where, in accordance with (ii) and (iii) of
Theorem~\ref{thmmain}, $\phi= q_0 g$, with $g$ the representation of
a~local
basic function~$g^*$. The integrand can then be written as $(p_0 +
q_0\Pi)g$, whence $S_0(P,Q) = \E_{P} g + \E_{Q}(\Pi g) = \E_{P^*} g^*
+ \E_{Q^*}(\Pi^* g^*)$---which thus has an invariant form,
$S_0^*(P^*,Q^*)$, independently of the representation employed.

We next demonstrate the corresponding property for the boundary term~$S_b$.

\begin{theorem}
\label{thmpB} We have
%
\begin{equation}
\label{eqspB}
\sum_r p_r B_r
= \sum_k \overline p_k \overline{B}_k \alpha{\circ} .
\end{equation}
\end{theorem}

\begin{pf}
On substituting (\ref{eqdiffq}), using (\ref{eqtform}) and
Theorem~\ref{thmdtrans}, and rearranging, the statement of the theorem
becomes
\begin{eqnarray*}
&&\sum_r p_r \sum_{m\geq r+1} \sum_{k=r+1}^m (-1)^{k-1-r} D^{k-1-r}
a_{mk}\,\frac{\partial}{\partial\overline{q}_m}\\
&&\qquad=\sum_r p_r \sum_{m\geq r+1} \sum_{k=r}^{m-1} a_{kr} (-1)^{m-1-k}
(\alpha
D)^{m-1-k} \,\frac{\partial}{\partial\overline{q}_m}\alpha{\circ }.
\end{eqnarray*}
The theorem will thus be proved if we can show
%
\begin{equation}
\label{eqBind}
{H_{m}^{(r)}\dvtx\sum_{k=r+1}^m (-1)^{k-1-r} D^{k-1-r}
a_{mk}\circ}
=\sum_{k=r}^{m-1} a_{kr} (-1)^{m-1-k} (\alpha
D)^{m-1-k}\alpha\circ\hspace*{-35pt}
\end{equation}
for all $m\geq r+1$. We prove (\ref{eqBind}) by induction on $m$.
First, $H_{r+1}^{(r)}$ holds since the left-hand side reduces to
$a_{r+1,r+1}\circ$ and the right-hand side reduces to
$a_{rr}\alpha\circ$, and these are equal by (\ref{eqarecur}). Now
suppose $H_m^{(r)}$ holds. Then
%
\begin{eqnarray}
&&\sum_{k=r+1}^m (-1)^{k-1-r} D^{k-1-r}
a_{mk}D\alpha{\circ}\\
\label{eqBind1}
&&\qquad=\sum_{k=r}^{m-1} a_{kr} (-1)^{m-1-k}
(\alpha D)^{m-1-k}\alpha D\alpha{\circ}.
\end{eqnarray}
But $ a_{mk}D = (D a_{mk}\circ) - a'_{mk}\circ$, so that, on
applying (\ref{eqarecur}), the left-hand side of (\ref{eqBind1}) becomes
\[
a_{mr}\alpha{\circ}- \sum_{k=r+1}^{m+1} (-1)^{k-1-r}
D^{k-1-r}a_{m+1,k}{\circ}.
\]
The right-hand side of (\ref{eqBind1}) straightforwardly becomes
\[
a_{mr}\alpha{\circ}- \sum_{k=r}^{m} a_{kr}(-1)^{m-k} (\alpha
D)^{m-k}\alpha{\circ},
\]
and we have thus verified $H_{m+1}$.
\end{pf}

It follows from (\ref{eqspB}) that $\sum_r p_r B_r q_0 =
\sum_r \overline p_r \overline B_r \overline q_0$, which thus defines
an invariant operator. Let now $S^*$, with representations $S$,
$\overline S$, derive from the $0$-homogeneous basic function $g^*$,
with representations $g$, $\overline g$. On using (\ref{eqsb}), in
which $\phi= q_0 g$, we get
\[
S_b(\bp,\bldq ) = \sum_r p_r B_r q_0 g = \sum_r
\overline p_r \overline B_r
\overline q_0 \overline g = \overline S_b(\overline{\bp},\overline
{\bldq }).
\]
Hence by (\ref{eqpbcphi}) the boundary contributions $S_\pm(P,Q)$ will
be the same in all representations.
%
\begin{expl}[(Modified Hyv\"{a}rinen rule)]
\label{exmod}
Take $\mathcal{X} = (0,\infty)$, $\overline\mathcal{X} = \R$,
$\gamma(x) \equiv\ln x$ (so $\overline X = \ln X$). Then $\alpha(x)
\equiv x$ and we find $\overline q_0 = x q_0$, $\overline q_1 = x
q_0 + x^2 q_1$, $\overline q_2 = x q_0 + 3x^2 q_1 + x^3 q_2$.

Let the scoring rule in the $\overline\xi$-representation,
$\overline S$, be defined by the Hyv\"{a}rinen formula:
%
\begin{equation}
\label{eqhyvbar}
\overline S(\overline x,\overline Q) = \frac{\overline
q''(\overline x)}{\overline q(\overline x)} -
\half\biggl\{\frac{\overline q'(\overline x)}{\overline
q(\overline x)}\biggr\}^2.
\end{equation}
This derives from the function $\overline g = -\half(\overline
q_1/\overline q_0)^2$.

Reexpressed in the $\xi$-representation, we have
%
\begin{equation}
\label{eqhyvnobar}
S( x, Q) = x^2\biggl[\frac{q''( x)}
{ q( x)} - \half\biggl\{\frac{ q'( x)}{q( x)}\biggr\}^2\biggr]
+ 2x \frac{q'(x)}{q(x)} + \half,
\end{equation}
which itself derives from the $\xi$-reexpression of $\overline g$,
viz., $g = -\half(1 + x q_1/q_0)^2$. That is, it is generated by
$\phi= q_0 g = -\half q_0 -x q_1 - \half x^2 q_1^2/q_0$. The
simpler choice $\phi^* = -\half x^2 q_1^2/q_0$ is equivalent to
$\phi$, and thus generates an equivalent scoring rule, with the same
divergence function; in fact, it simply eliminates the final term
$+\half$ in (\ref{eqhyvnobar}). This form of the scoring rule
also appears in equation (28) of \citet{Hyvarinenext}.

For this scoring rule, a~class ${\cal P}^*$ of distributions for the
basic variable~$X^*$ will be valid if, for $P, Q\in{\cal P}^*$,
$\overline p_0\{(\overline p_1/\overline p_0) - (\overline
q_1/\overline q_0)\}\rightarrow0$ as $\overline x\rightarrow\pm
\infty$, where these expressions are based on the $\overline
\xi$-representation [in which $\overline S$ is given by~(\ref{eqhyvbar})].
Reexpressing this in the $\xi$-representation, we
want
%
\begin{equation}
\label{eqnewb}
x^2 p_0\{( p_1/ p_0) - ( q_1/
q_0)\}\rightarrow0 \qquad\mbox{as }x \rightarrow0 \mbox{ or }\infty.
\end{equation}
At the lower end-point $0$ of $ {\cal X}$, this condition is less
restrictive than the corresponding condition (\ref{eqhyvdb}) for the
regular Hyv\"{a}rinen scoring rule defined directly on $ {\cal
X}$---although it becomes more restrictive at $\infty$.

In particular, suppose we consider the family $\mathcal{E}$ of
exponential densities:
\[
q(x \cd\theta) = \theta e^{-\theta x} \qquad(x, \theta> 0).
\]
For $p, q \in\mathcal{E}$, condition (\ref{eqnewb}) is satisfied,
whereas (\ref{eqhyvdb}) is not. If we tried to apply the unmodified
Hyv\"arinen score (\ref{eqhyv1}) to estimate $\theta$ in this model,
we would obtain $S(x, Q_\theta) = \half\theta^2$, and
(\ref{eqestimating}) would then appear to yield the clearly
nonsensical estimate $\widehat\theta\equiv0$. This is due to
failure of the boundary conditions, so that the original
Hyv\"{a}rinen rule is not in fact proper in this case. The modified
rule (\ref{eqhyvnobar}) is proper for this family, and yields the
consistent estimator $2\sum_i X_i/\sum_i X_i^2$.
\end{expl}

\section{Discussion and further work}

In this paper we have investigated local scoring rules only for the
case that the sample space is an open interval on the real line. The
general ideas extend to the case that the sample space is a~simply-connected $d$-dimensional differentiable manifold. This raises
challenging new technical problems, but could deliver a~fundamentally
improved understanding and illuminate issues associated with boundary
problems.

\section*{Acknowledgments}

We are grateful to Valentina Mameli for identifying errors and
simplifying proofs in an earlier version of this work. In addition we
have benefited from helpful discussions with Werner Ehm, Tilmann
Gneiting and Peter Whittle.



\printaddresses

\end{document}